\documentclass[11pt]{article}
\usepackage[T1]{fontenc}
\usepackage[utf8]{inputenc}
\usepackage{authblk}
\usepackage{psfrag,amsmath}
\usepackage{amssymb}
\usepackage{mathrsfs}
\usepackage{amsfonts}
\usepackage{amsmath}
\usepackage{mathrsfs,color}
\usepackage{epstopdf}
\usepackage{graphicx}
\usepackage{mathrsfs,amscd,amssymb,amsthm,amsmath,bm,url}

\makeatletter

\renewcommand{\@seccntformat}[1]{{\csname the#1\endcsname}{\normalsize .}\hspace{.5em}}
\makeatother

\def \[{\begin{equation}}
\def \]{\end{equation}}

\def \rank{{\rm rank}}
\def \Tr{{\rm Tr}}
\newtheorem{thm}{Theorem}[section]

\newtheorem{claim}{Claim}
\newtheorem{defi}{Definition}

\newtheorem{lem}[thm]{Lemma}
\newtheorem{cor}{Corollary}
\newtheorem{ex}[thm]{Example}
\newtheorem{pb}[thm]{Problem}

\newtheorem{conj}[thm]{Conjecture}

\begin{document}

\setlength{\baselineskip}{15pt}
\begin{center}{\Large \bf An arithmetic criterion for graphs being determined by their generalized $A_\alpha$-spectrum\footnote{Financially supported  by the National Natural Science Foundation of
China (Grant No. 11671164) and the Graduate Education Innovation Grant from Central China Normal University (Grant No. 2019CXPY001)}}
\vspace{4mm}

{\large Shuchao Li\footnote{Corresponding author. \\
\hspace*{5mm}{\it Email addresses}: lscmath@mail.ccnu.edu.cn (S.C. Li), \ wtsun2018@sina.com (W.T. Sun).},\ \ Wanting Sun}\vspace{4mm}

{\small Faculty of Mathematics and Statistics,  Central China Normal
University, Wuhan 430079, P.R. China}
\end{center}

\noindent {\bf Abstract}:\ Let $G$ be a graph on $n$ vertices, its adjacency matrix and degree diagonal matrix are denoted by $A(G)$ and $D(G)$, respectively. In 2017,  Nikiforov \cite{0007} introduced the matrix $A_{\alpha}(G)=\alpha D(G)+(1-\alpha)A(G)$ for  $\alpha\in [0, 1].$ The $A_\alpha$-spectrum of a graph $G$ consists of all the eigenvalues (including the multiplicities) of $A_\alpha(G).$ A graph $G$ is said to be determined by the generalized $A_{\alpha}$-spectrum (or, DGA$_\alpha$S for short) if whenever $H$ is a graph such that $H$ and $G$ share the same $A_{\alpha}$-spectrum and so do their complements, then $H$ is isomorphic to $G$. In this paper, when $\alpha$ is rational, we present a simple arithmetic condition for a graph being DGA$_\alpha$S. More precisely, put $A_{c_\alpha}:={c_\alpha}A_\alpha(G),$ here ${c_\alpha}$ is the smallest positive integer such that $A_{c_\alpha}$ is an integral matrix. Let $\tilde{W}_{{\alpha}}(G)=\left[{\bf 1},\frac{A_{c_\alpha}{\bf 1}}{c_\alpha},\ldots, \frac{A_{c_\alpha}^{n-1}{\bf 1}}{c_\alpha}\right]$, where ${\bf 1}$ denotes the all-ones vector. We prove that if  $\frac{\det \tilde{W}_{{\alpha}}(G)}{2^{\lfloor\frac{n}{2}\rfloor}}$ is an odd and square-free integer and the rank of $\tilde{W}_{{\alpha}}(G)$ is full over $\mathbb{F}_p$  for each odd prime divisor $p$ of $c_\alpha$, then $G$ is DGA$_\alpha$S except for even $n$ and odd $c_\alpha\,(\geqslant 3)$. By our obtained results in this paper we may deduce the main results in \cite{0005} and \cite{0002}.

\vspace{2mm} \noindent{\it Keywords:}
Generalized $A_{\alpha}$-spectrum; Cospectral; $A_{\alpha}$-matrix; Walk matrix
\vspace{2mm}


\section{\normalsize Introduction}\setcounter{equation}{0}
Throughout this paper, we assume $G$ is a simple (i.e., a finite, undirected, loopless and without multiple edges) connected graph, whose vertex set is $V_G=\{v_1,v_2,\ldots,v_n\}$ and edge set is $E_G$. The \textit{order} of $G$ is the number $n=|V_G|$ of its vertices and its \textit{size} is the number $|E_G|$ of its edges. Denote by $\bar{G}$ the complement of $G$. Unless otherwise stated, we follow the traditional notation and terminology (see \cite{0009}).

Given a graph $G$, its \textit{adjacency matrix} $A(G)$ is an $n\times n$ $0$-$1$ matrix whose $(i,j)$-entry is $1$ if and only if $v_i$ is adjacent to $v_j$ in $G$. 
Let $D(G)={\rm diag}(d_1,\ldots, d_n)$ be the diagonal matrix of vertex degrees in a graph $G$. The matrix $Q(G)=D(G)+A(G)$ is called the \textit{signless Laplacian matrix} of $G$ (see \cite{RD}). 
To track the gradual change of $A(G)$ into $Q(G),$ Nikiforov \cite{0007} introduced the \textit{$A_{\alpha}$-matrix} of a graph $G$, which is a convex combination of $D(G)$ and $A(G)$, that is,
$$
  A_{\alpha}(G)=\alpha D(G)+(1-\alpha)A(G),\ \ \ 0\leqslant \alpha\leqslant 1.
$$
One attractive feature of $A_{\alpha}$-matrices is that we can determine many of their properties from those of adjacency matrices or signless Laplacian matrices. Notice that $A_{\alpha}(G)$ is real and symmetric. Hence its eigenvalues are real. 
For short, the $A_{\alpha}$-spectral radius of $G$ (i.e., the largest eigenvalue of $A_{\alpha}(G)$) is called the \textit{$A_{\alpha}$-index} of $G.$ Note that
$$
  A(G)=A_0(G),\ \ \ Q(G)=2A_{\frac{1}{2}}(G) \ \ \ \text{and}\ \ \ D(G)=A_1(G).
$$

Recently, more and more people focused on the $A_{\alpha}$-matrix of a graph.  Nikiforov et al. \cite{Nik} gave some bounds on the $A_{\alpha}$-index of graphs, and they determined the unique tree with maximum (resp. minimum) $A_{\alpha}$-index among $n$-vertex trees. Nikiforov et al. \cite{005} and Xue et al. \cite{Xue}, independently, gave three edge graft transformations on $A_\alpha$-index. As applications, Xue et al. \cite{Xue} determined the graphs with maximum (resp. minimum) $A_{\alpha}$-index among all connected graphs with given diameter (resp. clique number). For more advances on the $A_{\alpha}$-spectra, we refer the reader to \cite{Chen,Huang,Li,LI2019,LI2020,001,Ni,Wang1,Wang,Xu2020} and references cited in.


It is well-known that the spectrum of a graph $G$ consists of all the eigenvalues (including the multiplicities) of the corresponding matrix associated with $G.$ In the literature, one usually studies the adjacency spectrum, Laplacian spectrum, signless Laplacian spectrum and $A_\alpha$-spectrum of a graph $G$, which are denoted by ${\rm Spec}_A(G),\,{\rm Spec}_L(G),\,{\rm Spec}_Q(G)$ and ${\rm Spec}_{\alpha}(G),$ respectively. We say two graphs are \textit{cospectral} if they share the same spectrum.

A graph $G$ is said to be \textit{determined by the spectrum} (DS for short) if, whenever $H$ is a graph such that $H$ and $G$ are cospectral, then $H$ is isomorphic to $G$ (here the matrix associated with $G$ should be clear in the context). In particular, for $\alpha\in[0,1),$ a graph $G$ is said to be \textit{determined by the generalized $A_\alpha$-spectrum} (or, DGA$_\alpha$S for short) if whenever $H$ is a graph such that ${\rm Spec}_{\alpha}(G) = {\rm Spec}_{\alpha}(H)$ and ${\rm Spec}_{\alpha}(\bar{G}) = {\rm Spec}_{\alpha}(\bar{H}),$ then $H$ is isomorphic to $G.$

``Which kinds of graphs are DS?" is a classical problem in spectral graph theory. The problem originates from chemistry and goes back to more than 60 years ago. In 1956, G\"{u}nthard and Primas \cite{9} raised the question in a paper that relates the theory of graph spectra to H\"{u}ckel's theory from chemistry. Kac \cite{11} asked a similar question: ``Can one hear the shape of a drum?". Fisher \cite{10} used the graph to model the shape of a drum. Then the sound of the drum can be identified by the eigenvalues of the corresponding graph. However, it turns out that determining whether a graph is DS is usually a difficult problem. We refer the reader to van Dam and Haemers \cite{van,van2} for some background and known results.


In the literature, many researchers studied the above problem in the context of the generalized $A_0$-spectrum and generalized $A_{\frac{1}{2}}$-spectrum (i.e., generalized adjacency spectrum and generalized $Q$-spectrum). Liu, Siemons and Wang \cite{Liu6} have constructed infinite families of graphs that are DG$A_0$S. Mao, Liu and Wang \cite{Mao3} gave a simple way to construct large DG$A_0$S graphs from small ones. Wang and Mao \cite{Wang4} presented a simple sufficient condition, under which they showed that $G\cup H$ is DG$A_0$S if and only if both $G$ and $H$ are DG$A_0$S. We make no attempt here to survey more important early contributions but instead refer the reader to \cite{A2,Liu5,0005,0001,0002,Wang01,0003,0004,Wang7} and references cited in.


Assume that $G$ is a graph with $A_\alpha$-matrix $A_\alpha(G)$ for $\alpha\in [0,1).$ In our whole context, we only consider that $\alpha$ is rational. Then let $c_\alpha$ be the smallest positive integer such that $A_{c_\alpha}(G):={c_\alpha} A_\alpha(G)$ is an integral matrix, that is, $c_\alpha$ is the smallest integer such that ${c_\alpha}\alpha$ and ${c_\alpha}(1-\alpha)$ are nonnegative integers. Define $W_{{\alpha}}(G)=[{\bf 1},A_{c_\alpha}{\bf 1},\ldots,A_{c_\alpha}^{n-1}{\bf 1}]$ to be the $A_{\alpha}$-walk matrix of $G$ for $\alpha\in[0,1)$, where ${\bf 1}$ denotes the all-ones vector and we always abbreviate $A_{c_\alpha}(G)$ to $A_{c_\alpha}$.

It is easy to see that $2^{\lfloor\frac{n}{2}\rfloor}c_\alpha^{n-1}$ is of course a factor of $\det W_{{\alpha}}(G)$ (based on Lemma \ref{lem3.4} below).
Furthermore, we shall notice that the rank of the $A_{\alpha}$-walk matrix $W_{{\alpha}}(G)$ is $1$ over $\mathbb{F}_{c_\alpha}$ if $c_\alpha>1.$ So, for $\alpha\in[0,1)$, we define the modified $A_{\alpha}$-walk matrix of $G$, written as $\tilde{W}_{{\alpha}}(G)$, to be $[{\bf 1},\frac{A_{c_\alpha}{\bf 1}}{c_\alpha},\ldots, \frac{A_{c_\alpha}^{n-1}{\bf 1}}{c_\alpha}].$ Recall that $A_{c_\alpha}(G){\bf 1} ={c_\alpha}A_\alpha(G){\bf 1}.$ It follows that $\tilde{W}_{{\alpha}}(G)$ is an integral matrix. 

Recently, Wang et al. \cite{0001,0002,0003,0004} and Qiu et al. \cite{0005} gave a simple arithmetic condition for a graph being DG$A_0$S and DG$A_{\frac{1}{2}}$S, respectively.
In this paper, we generalize their results to $A_\alpha$-matrix for $0\leqslant \alpha<1$. Our main result is given as follows. 
\begin{thm}\label{thm1.1}
Let $G$ be a graph with order $n\,(n\geqslant 5)$ and let $\alpha\in[0,1).$ If $\frac{\det \tilde{W}_{{\alpha}}(G)}{2^{\lfloor\frac{n}{2}\rfloor}}$ is odd and square-free, and the rank of $\tilde{W}_{{\alpha}}(G)$ is full over $\mathbb{F}_p$ for each odd prime $p\mid {c_\alpha},$ then $G$ is DGA$_{\alpha}$S except for even $n$ and odd $c_\alpha\,(\geqslant 3)$.
\end{thm}

The main idea of the proof of Theorem \ref{thm1.1} follows from Qiu \cite{0005} and Wang \cite{0001,0002}. Together with some new idea we make the proof work.

The remainder of the paper is organized as follows: In Section 2, we give some preliminary results that will be needed in the sequel. In Section 3, we present the proof of Theorem \ref{thm1.1}. In Section 4, we give examples of DGA$_\alpha$S graphs for some $\alpha\in[0,1)$. In Section 5, we give some conclusion remarks and some further research problems.
\section{\normalsize Preliminaries}\setcounter{equation}{0}
In this section, we present some preliminary results which will be needed in the subsequent sections.

If we focus on an integral matrix, its Smith Normal Form (SNF for short) is a very useful tool in our study. An integral matrix $V$ is said to be \textit{unimodular} if $\det V =\pm1.$ We firstly introduce the following well-known theorem.
\begin{thm}[\cite{0006}]\label{thm2.3}
Let $M$ be an integral matrix of order $n$. Then there exist unimodular matrices $V_1$ and $V_2$ such that $M = V_1SV_2$, where $S=\text{\rm diag}(s_1, s_2,\ldots,s_n)$ is the SNF of $M$ with $s_i\mid s_{i+1}$ for all $i\in \{1,2,\ldots,n-1\},$ and $s_i$ is called the $i$-th elementary divisor of $M$.
\end{thm}
Notice that the SNF of a matrix can be computed efficiently (see \cite[P50]{0008}). The following lemma is the key result in the proof of our main result.
\begin{lem}[\cite{0001}]\label{lem2.4}
With the above notations, the system of congruence equations $M{\bf x}\equiv {\bf 0} \pmod{p^2}$ has a solution ${\bf x}\not\equiv{\bf 0} \pmod{p}$ if and only if $p^2 \mid s_n.$
\end{lem}

Next, we will illustrate the main strategy to prove a graph to be DGA$_\alpha$S for $\alpha\in[0,1).$ 
The following lemma is the groundwork of our method, which gives a simple characterization of two
graphs having the same generalized $A_\alpha$-spectrum. It is an analogous result of adjacency matrix in \cite{0004} and signless Laplacian matrix in \cite{0005}. For convenience, we always set $a:=c_\alpha\alpha$ and $b:=c_\alpha(1-\alpha)$ and denote by $O_{n}(\mathbb{Q})$ the set of all $n\times n$ rational orthogonal matrices.

\begin{lem}\label{thm2.1}
Let $G$ be a graph of order $n$ such that $\det W_{{\alpha}}(G)\neq0$ for $\alpha\in[0,1).$ Then there exists a graph $H$ such that $G$ and $H$ share the same generalized $A_{\alpha}$-spectrum if and only if there exists a unique matrix $U\in O_{n}(\mathbb{Q})$ satisfying that
\[\label{eq:2.1}
  U^TA_{c_\alpha}(G)U=A_{c_\alpha}(H)\ \  \text{and}\ \ U{\bf 1}={\bf 1}.
\]
\end{lem}
\begin{proof}
\textit{Sufficiency.}\ Assume that there exists a matrix $U\in O_{n}(\mathbb{Q})$ satisfying \eqref{eq:2.1}. It is routine to check that
\[\label{eq:2.01}
  A_{c_\alpha}(\bar{G})=aD(\bar{G})+bA(\bar{G})=bJ+(a(n-1)-b)I-A_{c_\alpha}(G),
\]
where $I$ and $J$ denote the identity matrix and the all one's matrix, respectively. Notice that $U{\bf 1}={\bf 1}$ and $U\in O_{n}(\mathbb{Q}).$ Therefore,
$$
  U^TA_{c_\alpha}(\bar{G})U=U^T(bJ+(a(n-1)-b)I-A_{c_\alpha}(G))U=bJ+(a(n-1)-b)I-A_{c_\alpha}(H)=A_{c_\alpha}(\bar{H}).
$$
Hence $G$ and $H$ are cospectral with respect to the generalized $A_{\alpha}$-spectrum.

\textit{Necessity.}\
Note that ${\rm Spec}_\alpha(G)={\rm Spec}_\alpha(H)$ and ${\rm Spec}_\alpha(\bar{G})={\rm Spec}_\alpha(\bar{H}).$ 
Together with \eqref{eq:2.01}, one obtains that
\[\label{eq:03}
\det(\lambda I-A_{c_\alpha}(G))=\det(\lambda I-A_{c_\alpha}(H))\ \text{and}\ \det(\lambda I+bJ-A_{c_\alpha}(G))=\det(\lambda I+bJ-A_{c_\alpha}(H))
\]
hold for all real $\lambda.$
Furthermore, it is routine to check that
\begin{align}\notag
  \det(\lambda I+bJ-A_{c_\alpha}(G))&=\det(\lambda I-A_{c_\alpha}(G)+b{\bf 1}{\bf 1}^T)\\\notag
  &=\det(\lambda I-A_{c_\alpha}(G))\det(I+b(\lambda I-A_{c_\alpha}(G))^{-1}{\bf 1}{\bf 1}^T)\\\label{eq:01}
  &=(1+b{\bf 1}^T(\lambda I-A_{c_\alpha}(G))^{-1}{\bf 1})\det(\lambda I-A_{c_\alpha}(G))
\end{align}
for all $\lambda\not\in \sigma(A_{c_\alpha}(G)),$ where $\sigma(M)$ denotes the set of all distinct eigenvalues of the matrix $M.$ Similarly, for all $\lambda\not\in \sigma(A_{c_\alpha}(H)),$ we have
\[\label{eq:02}
  \det(\lambda I+bJ-A_{c_\alpha}(G))=(1+b{\bf 1}^T(\lambda I-A_{c_\alpha}(H))^{-1}{\bf 1})\det(\lambda I-A_{c_\alpha}(H)).
\]
Combining \eqref{eq:03}-\eqref{eq:02}, we have
\[\label{eq:2.2}
  {\bf 1}^T(\lambda I-A_{c_\alpha}(G))^{-1}{\bf 1}={\bf 1}^T(\lambda I-A_{c_\alpha}(H))^{-1}{\bf 1}.
\]
Notice that $A_{c_\alpha}(G)$ is a real symmetric matrix. Hence all the linearly independent eigenvectors of $A_{c_\alpha}(G)$ form an orthonormal basis of $\mathbb{R}_n.$ Group the eigenvectors with respect to the eigenvalue $\mu$ into matrix $P_{\mu}$ for each $\mu\in \sigma({A_{c_\alpha}}(G))$. Without loss of generality, assume that $\sigma({A_{c_\alpha}}(G))=\sigma(A_{c_\alpha}(H))=\{\mu_1,\mu_2,\ldots,\mu_s\}.$ Therefore,
\begin{equation*}
    A_{c_\alpha}(G)\left[P_{\mu_1}, P_{\mu_2}, \ldots, P_{\mu_s}\right]=\left[P_{\mu_1}, P_{\mu_2}, \ldots, P_{\mu_s}\right]\left[
                       \begin{array}{cccc}
                         \mu_1I_1 &  &  &  \\
                          & \mu_2I_2 &  &  \\
                          &  & \ddots &  \\
                          &  &  & \mu_sI_s \\
                       \end{array}
                     \right],
\end{equation*}
where $I_i$ denotes the identity matrix whose order equals the multiplicity of $\mu_i$ for $1\leqslant i\leqslant s.$ Hence, for each $\lambda\not\in \sigma(A_{c_\alpha}(G)),$ we obtain
$$
  (\lambda I-A_{c_\alpha}(G))^{-1}=\sum_{i=1}^s\frac{1}{\lambda-\mu_i}P_{\mu_i} P_{\mu_i}^T.
$$
By a similar discussion, one has, for each $\lambda\not\in \sigma(A_{c_\alpha}(H)),$ that
$$
 (\lambda I-A_{c_\alpha}(H))^{-1}=\sum_{i=1}^s\frac{1}{\lambda-\mu_i}R_{\mu_i} R_{\mu_i}^T,
$$
here $R_{\mu_i}(1\leqslant i\leqslant s)$ denotes the matrix for $H$ corresponding to $P_{\mu_i}$ for $G.$
In view of \eqref{eq:2.2}, we obtain, for each $\lambda\not\in \sigma(A_{c_\alpha}(G)),$ that
$$
  \sum_{i=1}^s\frac{\left\|P_{\mu_i}^T{\bf 1}\right\|}{\lambda-\mu_i}=\sum_{i=1}^s\frac{\left\|R_{\mu_i}^T{\bf 1}\right\|}{\lambda-\mu_i}.
$$
This implies that $\left\|P_{\mu_i}^T{\bf 1}\right\|=\left\|R_{\mu_i}^T{\bf 1}\right\|$ for all $i\in\{1,\ldots,s\}.$ Hence, there exists an orthogonal matrix $H_{\mu_i}$ such that $P_{\mu_i}^T{\bf 1}=H_{\mu_i} R_{\mu_i}^T{\bf 1}$ for all $i\in\{1,\ldots,s\}.$ Let
$$
  U=[P_{\mu_1},P_{\mu_2},\ldots,P_{\mu_s}][ R_{\mu_1}H_{\mu_1}^T, R_{\mu_2}H_{\mu_2}^T,\ldots, R_{\mu_s}H_{\mu_s}^T]^T.
$$
It is straightforward to check that $U$ is an orthogonal matrix such that $U{\bf 1}={\bf 1}$ and $U^TA_{c_\alpha}(G)U=A_{c_\alpha}(H).$ Therefore, $U^TA_{c_\alpha}^k(G){\bf 1}=A_{c_\alpha}^k(H){\bf 1}$ for each $k\in\{0,1,\ldots,n-1\},$ which yields that $U^TW_{{\alpha}}(G)=W_{{\alpha}}(H).$ It follows from  $\det W_{\alpha}(G)\neq0$ that $\det W_{\alpha}(H)\neq0,$ and thus $U=W_{{\alpha}}(G)W_{{\alpha}}(H)^{-1}$ is a rational orthogonal matrix satisfying \eqref{eq:2.1}.

Now, we show the uniqueness of $U.$ Suppose to the contrary that there exist two distinct matrices $U_1,\,U_2\in O_n(\mathbb{Q})$ satisfying \eqref{eq:2.1}. Then we obtain $U_1^TW_{{\alpha}}(G)=U_2^TW_{{\alpha}}(G)=W_{{\alpha}}(H).$ Notice that $\det W_{{\alpha}}(G)\neq 0.$ Hence $U_1=U_2,$ a contradiction.

 This completes the proof.
\end{proof}
For $\alpha\in[0,1),$ we now define the following notation that will be used frequently in the sequel:
$$
  \Gamma_\alpha(G)=\{U:U\in O_n(\mathbb{Q}),\,U^T{A_{c_\alpha}}(G)U ={A_{c_\alpha}}(H)\,\text{for some graph}\,H\, \text{and}\, U{\bf 1} = {\bf 1}\}.
$$

The next lemma extends the results for the adjacency matrix and signless Laplacian matrix found in standard texts like \cite{0005,0004} to $A_\alpha$-matrix.
\begin{lem}\label{thm2.2}
Let $G$ be a graph and let $\alpha\in[0,1)$. If $\det W_{{\alpha}}(G)\neq 0,$ then $G$ is DGA$_{\alpha}$S if and only if each matrix in $\Gamma_\alpha(G)$ is a permutation matrix.
\end{lem}
\begin{proof}
\textit{Sufficiency.}\ Suppose to the contrary that $G$ is not DGA$_{\alpha}$S. Then there exists a graph $H\not\cong G$ such that $G$ and $H$ share the same generalized $A_\alpha$-spectrum. 
According to Lemma \ref{thm2.1}, there is a unique matrix $U\in O_n(\mathbb{Q})$ satisfying \eqref{eq:2.1}. Therefore, $U\in \Gamma_\alpha(G).$ Notice that $H\not\cong G.$ Hence $U$ is not a permutation matrix, which contradicts the assumption that each matrix in $\Gamma_\alpha(G)$ is a permutation matrix.

\textit{Necessity.}\ Suppose that there exists a matrix $U\in \Gamma_\alpha(G)$ and $U$ is not a permutation matrix. Assume  that $H$ is a graph satisfying  $U^T{A_{c_\alpha}}(G)U ={A_{c_\alpha}}(H).$ Based on the uniqueness of $U,$ we obtain that $H\not\cong G,$ a contradiction.
\end{proof}
Further on we need the following definition.
\begin{defi}\label{defi1}
  Let $U$ be a rational orthogonal matrix. The \textit{level} of $U$, denoted by $l(U)$ (or simply by $l$ if there is no danger of ambiguity), is the smallest positive integer $k$ such that $kU$ is an integral matrix.
\end{defi}

Obviously, $l(U)$ is the least common denominator of all the entries of the matrix $U$. It is routine to check that a matrix $U\in O_{n}(\mathbb{Q})$ with $U{\bf 1}={\bf 1}$ is a permutation matrix if and only if $l(U)=1.$ In view of Lemma~\ref{thm2.2}, for a given graph $G$ and $\alpha\in[0,1),$ our main strategy in proving that $G$ is DGA$_{\alpha}$S is equivalent to show that $l(U)=1$ for all $U\in \Gamma_\alpha(G).$


In the remaining part of this section, we give a technical lemma, which also plays an important role in the proof of Theorem \ref{thm1.1}. In what follows, when there is no scope for ambiguity, we always put $A:=A(G),$ $D:=D(G),$ $A_{c_\alpha}:=A_{c_\alpha}(G)$ and  $\tilde{W}_\alpha:=\tilde{W}_\alpha(G).$ 
\begin{lem}\label{lem2.5}
Let $G$ be a graph on $n$ vertices and let $\alpha\in[0,1)$. Then ${\bf 1}^TA_{c_\alpha}{\bf 1}\equiv 0 \pmod{2{c_\alpha}}$ and ${\bf 1}^TA_{c_\alpha}^k{\bf 1}\equiv 0 \pmod{2{c_\alpha^2}}$ for each integer $k\geqslant 2.$
\end{lem}
\begin{proof}
Notice that $D{\bf 1}=A{\bf 1}={\bf d},$ here ${\bf d}=(d_1,d_2,\ldots, d_n)^T$ and $d_i$ denotes the degree of the $i$-th vertex of $G.$ Hence, $A_{c_\alpha}{\bf 1}=(a A+bD){\bf 1}={c_\alpha}{\bf d}.$ Therefore,
$$
  {\bf 1}^TA_{c_\alpha}{\bf 1}={c_\alpha}{\bf 1}^T{\bf d}={c_\alpha}\sum_{i=1}^nd_i=2{c_\alpha}|E_G|\equiv 0 \pmod{2{c_\alpha}}.
$$
Furthermore,
$$
  {\bf 1}^TA_{c_\alpha}^2{\bf 1}=(A_{c_\alpha}{\bf 1})^T(A_{c_\alpha}{\bf 1})={c_\alpha^2}\sum_{i=1}^nd_i^2\equiv{c_\alpha^2}\sum_{i=1}^nd_i=2{c_\alpha^2}|E_G|\equiv0\pmod{2{c_\alpha^2}}.
$$

In what follows, we prove that ${\bf 1}^TA_{c_\alpha}^k{\bf 1}\equiv 0 \pmod{2{c_\alpha^2}}$ for each integer $k\geqslant 3.$ Recall that $A_{c_\alpha}{\bf 1}={c_\alpha}{\bf d}={c_\alpha}D{\bf 1}.$ Therefore, 
$$
  {\bf 1}^TA_{c_\alpha}^k{\bf 1}={c_\alpha^2}{\bf 1}^TD(aA+bD)^{k-2}D{\bf 1}\equiv{c_\alpha^2}\Tr(D(aA+bD)^{k-2}D)\pmod{2{c_\alpha^2}},
$$
here $\Tr(D(aA+bD)^{k-2}D)$ denotes the trace of $D(aA+bD)^{k-2}D).$ Hence, it suffices to show that $\Tr(D(aA+bD)^{k-2}D)\equiv0\pmod{2}$ for each integer $k\geqslant 3.$

For the ease of expression, we define $\mathfrak{X}$ as the free monoid generated by $\{x, y\},$ and
$$
  \mathfrak{X}_m=\{X\in \mathfrak{X}:\text{the length of}\, X\, \text{is}\, m\}.
$$
Define a mapping $\tau$ on $\mathfrak{X}_m$ such that $\tau(X)=M_mM_{m-1}\ldots M_1:=X^{\tau},$ where $M_i$ is the $i$-th character of $X$ for each $i\in \{1, 2,\ldots,m\}.$ Denote by $\underline{X}=M_1M_2\cdots M_m$ the product of the string of matrices in $X$, where $M_i=aA$ if the $i$-th character of $X$ is $x$, and $M_i=bD$ if the $i$-th character of $X$ is $y$, for $1\leqslant i\leqslant m.$ It is routine to check that $\underline{X}^T=\underline{X^{\tau}}$ and $X^{\tau}\in \mathfrak{X}_m$ is uniquely determined by $X.$

Using the above notations, one has
$$
  \Tr(D(aA+bD)^{k-2}D)=\sum_{{X}\in \mathfrak{X}_{k-2}}\Tr(D\underline{X}D).
$$
Notice that $\Tr(D\underline{X}D)=\Tr((D\underline{X}D)^T)=\Tr(D\underline{X^{\tau}}D).$ Therefore,
$$
  \sum_{{X}\in \mathfrak{X}_{k-2},\,X\neq X^{\tau}}\Tr(D\underline{X}D)\equiv 0\pmod{2}.
$$
Hence,
$$
  \Tr(D(aA+bD)^{k-2}D)\equiv\sum_{{X}\in \mathfrak{X}_{k-2},\,X=X^{\tau}}\Tr(D\underline{X}D)\pmod{2}.
$$
Now, we proceed by distinguishing the parity on $k$.

{\bf Case 1.}\ $k$ is even. Recall that $A{\bf 1}=D{\bf 1}.$ Hence,
\begin{align*}
  \Tr(D\underline{X^{\tau}}AA\underline{X}D)&=\Tr(A\underline{X}DD\underline{X^{\tau}}A)\equiv {\bf 1}^TA\underline{X}DD\underline{X^{\tau}}A{\bf 1}\\
  &={\bf 1}^TD\underline{X}DD\underline{X^{\tau}}D{\bf 1}\equiv \Tr(D\underline{X}DD\underline{X^{\tau}}D)\pmod{2}.
\end{align*}
It follows that
\begin{align*}
    \Tr(D(aA+bD)^{k-2}D)&\equiv\sum_{{X}\in \mathfrak{X}_{k-2},\,X=X^{\tau}}\Tr(D\underline{X}D)=\sum_{X\in \mathfrak{X}_{\frac{k}{2}-2}}\Tr(D\underline{X^{\tau}}(aAaA+bDbD)\underline{X}D)\\
    &=a^2\sum_{X\in \mathfrak{X}_{\frac{k}{2}-2}}\Tr(D\underline{X^{\tau}}AA\underline{X}D)+b^2\sum_{X\in \mathfrak{X}_{\frac{k}{2}-2}}\Tr(D\underline{X^{\tau}}DD\underline{X}D)\\
    &=(a^2+b^2)\sum_{X\in \mathfrak{X}_{\frac{k}{2}-2}}\Tr(D\underline{X^{\tau}}DD\underline{X}D)\\
    &\equiv(a+b)\sum_{X\in \mathfrak{X}_{\frac{k}{2}-2}}\Tr(D\underline{X^{\tau}}DD\underline{X}D)\pmod{2}.
\end{align*}
The last congruence equation follows from the fact that $t^2\equiv t\pmod{2}$ for each positive integer $t.$

If $a+b$ is even, then $\Tr(D(aA+bD)^{k-2}D)\equiv0\pmod{2},$ as desired. If $a$ is even and $b$ is odd, then
$$
  \Tr(D(aA+bD)^{k-2}D)\equiv \Tr(b^{k-2}D^k)=b^{k-2}\sum_{i=1}^n d_i^k\equiv b^{k-2}\sum_{i=1}^n d_i=2b^{k-2}|E_G|\equiv0\pmod{2},
$$
as desired. If $a$ is odd and $b$ is even, then
\begin{align*}
  \Tr(D(aA+bD)^{k-2}D)&\equiv a^{k-2}\Tr(DA^{k-2}D)\equiv {a{\bf 1}^TDA^{k-2}D{\bf 1}}=a{\bf 1}^TA^{k}{\bf 1}\\
&={a\sum_{i,j}(A^{k-1})_{ij}(A)_{ij}}= 2a\sum_{1\leqslant i<j\leqslant n}(A^{k-1})_{ij}(A)_{ij}\equiv 0\pmod {2},
\end{align*}
as desired.

{\bf Case 2.}\ $k$ is odd. It is straightforward to check that
\begin{align*}
    \Tr(D\underline{X^{\tau}}A\underline{X}D)&=\Tr(\underline{X}D D\underline{X^{\tau}}A)=\sum_{i,j}(\underline{X}D D\underline{X^{\tau}})_{ij}(A)_{ij}\\
    &=2\sum_{1\leqslant i<j\leqslant n}(\underline{X}D D\underline{X^{\tau}})_{ij}(A)_{ij}\equiv0\pmod{2}.
\end{align*}
Thus,
\begin{align*}
    \Tr(D(aA+bD)^{k-2}D)&\equiv\sum_{{X}\in \mathfrak{X}_{k-2},\,X=X^{\tau}}\Tr(D\underline{X}D)=\sum_{X\in \mathfrak{X}_{\frac{k-3}{2}}}\Tr(D\underline{X^{\tau}}(aA+bD)\underline{X}D)\\
    &\equiv b\sum_{X\in \mathfrak{X}_{\frac{k-3}{2}}}\Tr(D\underline{X^{\tau}}D\underline{X}D)\pmod{2}.
\end{align*}
It is easy to see that
$$
  \Tr(D\underline{X^{\tau}}D\underline{X}D)\equiv\Tr(D\underline{X^{\tau}}DD\underline{X}D)=\Tr(D\underline{X}DD\underline{X^{\tau}}D)
  \equiv\Tr(D\underline{X}D\underline{X^{\tau}}D)\pmod{2}.
$$
Then, it follows that
\begin{align*}
    \Tr(D(aA+bD)^{k-2}D)&\equiv b\sum_{X\in \mathfrak{X}_{\frac{k-3}{2}}}\Tr(D\underline{X^{\tau}}D\underline{X}D)\equiv b\sum_{X\in \mathfrak{X}_{\frac{k-3}{2}},\,X=X^{\tau}}\Tr(D\underline{X^{\tau}}D\underline{X}D)\\
    &\equiv b\sum_{X\in \mathfrak{X}_{\frac{k-3}{2}},\,X=X^{\tau}}\Tr(D\underline{X^{\tau}}DD\underline{X}D)\equiv b\sum_{X\in \mathfrak{X}_{\frac{k-3}{2}},\,X=X^{\tau}}({\bf 1}^TD\underline{X^{\tau}}D)(D\underline{X}D{\bf 1})\\
    &\equiv b\sum_{X\in \mathfrak{X}_{\frac{k-3}{2}},\,X=X^{\tau}}{\bf 1}^TD\underline{X}D{\bf 1}
    \equiv b\Tr(D(aA+bD)^{\frac{k-3}{2}}D)\pmod{2}.
\end{align*}
Then by induction on $k$, we can show that $\Tr(D(aA+bD)^{k-2}D)\equiv0\pmod{2}.$

Combining Cases 1-2, we obtain that ${\bf 1}^TA_{c_\alpha}^k{\bf 1}\equiv 0 \pmod{2{c_\alpha^2}}$ for each integer $k\geqslant 3.$ This completes the proof.
\end{proof}
\section{\normalsize Proof of Theorem \ref{thm1.1}}\setcounter{equation}{0}
In this section, we give the proof of Theorem \ref{thm1.1}.  In view of Lemma \ref{thm2.2}, it suffices to prove that for each $U\in \Gamma_\alpha(G)$, the conditions of Theorem \ref{thm1.1} mean  that the level $l(U)=1,$ which is equivalent to show that any prime $p$ is not a divisor of $l(U).$

In what follows, we will use the finite field notation $\mathbb{F}_p$ and mod $p$ (for a prime $p$) interchangeably, and denote by $\rank_p(M)$ the \textit{rank} of an integral matrix $M$ over $\mathbb{F}_p.$

Let $\mathfrak{F}_n$ be the set of all graphs $G$ on $n$ vertices satisfying that $\frac{\det \tilde{W}_{{\alpha}}(G)}{2^{\lfloor\frac{n}{2}\rfloor}}$ is odd and square-free, and $\rank_p(\tilde{W}_{{\alpha}}(G))=n$ for each odd prime $p\mid{c_\alpha}$ and $\alpha\in[0,1)$

In order to complete the proof of Theorem \ref{thm1.1}, we need the following key results.
\begin{thm}\label{thm3.1}
Let $G\in \mathfrak{F}_n$ and let $\alpha\in[0,1).$ If $U\in \Gamma_\alpha(G)$ with level $l$ and $p$ is an odd prime, then $p\nmid l.$
\end{thm}
\begin{thm}\label{thm3.3}
Let $G\in \mathfrak{F}_n\,(n\geqslant 5)$ and let $\alpha\in[0,1).$ If $U\in \Gamma_\alpha(G)$ with level $l$, then $l$ is odd except for even $n$ and odd $c_\alpha\,(\geqslant 3)$.
\end{thm}
We postpone the proofs of Theorems \ref{thm3.1} and \ref{thm3.3} to the subsequent of this section. As a corollary, we present the proof of Theorem \ref{thm1.1}:
\begin{proof}[\bf Proof of Theorem \ref{thm1.1}]\
Combining Theorems \ref{thm3.1} and \ref{thm3.3}, our results follows immediately.
\end{proof}
In Theorem 1.1, put $\alpha=0$, then $c_\alpha=1$ and we may obtain the following corollary immediately. It gives a  simple arithmetic condition for determining whether a graph is determined by the generalized adjacency spectrum, which  was obtained by Wang \cite{0002}.
\begin{cor}\label{cor1.01}
Let $G$ be a graph of order $n$. If $\det \frac{\tilde{W}_0(G)}{2^{\lfloor\frac{n}{2}\rfloor}}$ is odd and square-free, then $G$ is DG$A_0$S.
\end{cor}
In Theorem 1.1, put $\alpha=\frac{1}{2}$, then $c_\alpha=2$ and we may obtain the following corollary immediately. It gives a  simple arithmetic condition for determining whether a graph is determined by the generalized $Q$-spectrum, which was obtained by Qiu, Ji and Wang \cite{0005}.
\begin{cor}\label{cor1.02}
Let $G$ be a graph of order $n$. If $\det \frac{\tilde{W}_\frac{1}{2}(G)}{2^{\lfloor\frac{n}{2}\rfloor}}$ is odd and square-free, then $G$ is DG$A_\frac{1}{2}$S.
\end{cor}

\subsection{\normalsize Proof of Theorem \ref{thm3.1}}
In this subsection, we are devoted to the proof of Theorem \ref{thm3.1}. 
Before giving the proof of Theorem \ref{thm3.1}, we present the following needed lemmas.
\begin{lem}\label{lem3.2}
Let $G$ be a graph and let $\alpha\in[0,1).$ Assume that $U\in \Gamma_\alpha(G)$ with level $l$ and $p$ is a prime division of $l.$ Then for each integer $k\geqslant 0,$ there exists an integral column vector ${\bf v}\not\equiv{\bf 0}\pmod{p}$ satisfying that
\[\label{eq:3.2}
  {\bf v}^TA_{c_\alpha}^k{\bf v}\equiv 0\pmod{p^2}\ \ \text{and}\ \ \tilde{W}_{{\alpha}}^T{\bf v}\equiv{\bf 0}\pmod{p}.
\]
\end{lem}
\begin{proof}
Let $H$ be a graph such that $U^TA_{c_\alpha}(G)U=A_{c_\alpha}(H)$ and let $\bar{U}=lU.$ Based on Definition \ref{defi1}, there exists a column ${\bf v}$ of $\bar{U}$ such that ${\bf v}\not\equiv{\bf 0}\pmod{p}.$ It is easy to check that $\bar{U}^T{A_{c_\alpha}^k}(G)\bar{U}=l^2{A_{c_\alpha}^k}(H)\equiv{\bf 0}\pmod{p^2}.$ It follows that ${\bf v}^TA_{c_\alpha}^k(G){\bf v}\equiv 0\pmod{p^2}$ for each integer $k\geqslant 0.$

Notice that $U{\bf 1}={\bf 1}.$ Hence, $U^TA_{c_\alpha}^k(G){\bf 1}=A_{c_\alpha}^k(H){\bf 1}$ for each integer $k\geqslant 0.$ It follows that  $U^T\frac{A_{c_\alpha}^k(G){\bf 1}}{c_\alpha}=\frac{A_{c_\alpha}^k(H){\bf 1}}{c_\alpha}$ for each integer $k\geqslant 1,$ and therefore $U^T\tilde{W}_{{\alpha}}(G)=\tilde{W}_{{\alpha}}(H)$ is an integral matrix. Hence, $\tilde{W}_{{\alpha}}(G)^T{\bf v}\equiv{\bf 0}\pmod{p}.$ This completes the proof.
\end{proof}
\begin{lem}\label{lem3.6}
Let $G$ be a graph with $\det \tilde{W}_{{\alpha}}(G)\neq 0$ and let $\alpha\in[0,1)$. If $U\in \Gamma_\alpha(G)$ with level $l,$ then $l\mid s_n,$ where $s_n$ is the $n$-th elementary division of $\tilde{W}_{{\alpha}}(G).$
\end{lem}
\begin{proof}
Let $H$ be a graph such that $U^T{A_{c_\alpha}}(G)U={A_{c_\alpha}}(H).$ In view of the proof of Lemma~\ref{lem3.2}, one has 
$U^T\tilde{W}_{{\alpha}}(G)=\tilde{W}_{{\alpha}}(H).$ That is to say, $U^T=\tilde{W}_{{\alpha}}(H)\tilde{W}_{{\alpha}}(G)^{-1}.$

By Theorem \ref{thm2.3}, there exist unimodular matrices $V_1$ and $V_2$ such that $\tilde{W}_{{\alpha}}(G)=V_1SV_2,$ where $S={\rm diag}(s_1,s_2,\ldots,s_n)$ and $s_i\mid s_{i+1}$ for each $i\in\{1,\ldots,n-1\}.$ Thus,
$$
  s_nU^T=\tilde{W}_{{\alpha}}(H)V_2^{-1}{\rm diag}\left(\frac{s_n}{s_1},\frac{s_n}{s_2},\ldots,\frac{s_n}{s_n}\right)V_1^{-1}.
$$
Notice that $V_1^{-1}$ and $V_2^{-1}$ are integral matrices. Hence, $s_nU^T$ is an integral matrix. Therefore, $l\mid s_n.$ This completes the proof.
\end{proof}
Now, we are ready to prove Theorem \ref{thm3.1}.
\begin{proof}[\bf Proof of Theorem \ref{thm3.1}]\
Firstly, we consider $p\mid c_\alpha.$ It follows from Lemma \ref{lem3.6} that $p\nmid l.$ Otherwise, $p\mid \det \tilde{W}_\alpha$ and therefore $\rank_p(\tilde{W}_{{\alpha}})\neq n,$ which contradicts the definition of $\mathfrak{F}_n.$ So, in what follows, it suffices to assume that $p\nmid c_\alpha.$

Suppose, to the contrary, that $p\mid l.$ In view of Theorem \ref{thm2.3}, there exist unimodular matrices $V_1$ and $V_2$ such that $\tilde{W}_{{\alpha}}=V_1SV_2,$ where $S={\rm diag}(s_1,s_2,\ldots,s_n)$ and $s_i\mid s_{i+1}$ for each $i\in\{1,\ldots,n-1\}.$ By Lemma \ref{lem3.6}, one has $p\mid s_n.$

Notice that $\det \tilde{W}_{{\alpha}}=\pm \det S.$ Therefore $p\mid \det \tilde{W}_{{\alpha}}.$ Based on the definition of $\mathfrak{F}_n,$ one has $p^2\nmid \det \tilde{W}_{{\alpha}}.$ Hence, $\rank_p(\tilde{W}_{{\alpha}})=n-1.$ In view of the proof Lemma~\ref{lem3.2}, we know that each column of $\bar{U}$ satisfies \eqref{eq:3.2}. Thus, $\tilde{W}_{{\alpha}}^T\bar{U}\equiv{\bf 0}\pmod{p}$ and so $\rank_p(\bar{U})=1.$ It follows that there exists an integral vector ${\gamma}$ such that ${\bf v}{ \gamma}^T\equiv \bar{U}\pmod{p},$ where ${\bf v}$ is the $j$-th column of $\bar{U}$ satisfying both ${\bf v}\not\equiv{\bf 0}\pmod{p}$ and \eqref{eq:3.2}. Let $H$ be a graph such that $U^TA_{c_\alpha}(G)U=A_{c_\alpha}(H).$ Thus,
$$
  A_{c_\alpha}(G){\bf v}=\bar{U}A_{c_\alpha}(H)_j\equiv {\bf v}(\gamma^TA_{c_\alpha}(H)_j)=\lambda_0{\bf v}\pmod{p},
$$
here $A_{c_\alpha}(H)_j$ denotes the $j$-th column of $A_{c_\alpha}(H)$ and $\lambda_0=\gamma^TA_{c_\alpha}(H)_j.$ Hence, $\rank_p(A_{c_\alpha}(G)-\lambda_0I)\neq n.$ We proceed by considering the following three cases.

{\bf Case 1.}\ $\rank_p(A_{c_\alpha}-\lambda_0I)= n-1.$ Notice that ${\bf v}^T(A_{c_\alpha}-\lambda_0I)\equiv{\bf 0}\pmod{p},$ ${\bf v}^T{\bf v}\equiv0\pmod{p^2}$ and ${\bf v}^T{\bf 1}\equiv0\pmod{p}.$ Hence, there exist integral vectors ${\bf y}$ and ${\bf u}$ satisfying ${\bf v}\equiv(A_{c_\alpha}-\lambda_0I){\bf y}\pmod{p}$ and ${\bf 1}\equiv(A_{c_\alpha}-\lambda_0I){\bf u}\pmod{p}.$ That is to say, ${\bf 1}=(A_{c_\alpha}-\lambda_0I){\bf u}+p\beta$ for some integral vector $\beta.$ It follows that
\begin{align}\label{eq:3.3}
    \tilde{W}_{{\alpha}}=\left[{\bf 1},\frac{A_{c_\alpha}{\bf 1}}{c_\alpha},\ldots,\frac{A_{c_\alpha}^{n-1}{\bf 1}}{c_\alpha}\right]=(A_{c_\alpha}-\lambda_0I)X
    +p\left[\beta,\frac{A_{c_\alpha}\beta}{c_\alpha},\ldots,\frac{A_{c_\alpha}^{n-1}\beta}{c_\alpha}\right],
\end{align}
here $X=\left[{\bf u},\frac{A_{c_\alpha}{\bf u}}{c_\alpha},\ldots,\frac{A_{c_\alpha}^{n-1}{\bf u}}{c_\alpha}\right].$ Notice that $X$ does not need to be an integral matrix, but $X\pmod{p}$ is meaningful since $p$ is an odd prime and $p\nmid c_\alpha$. Therefore,
$$
  {c_\alpha}\frac{\tilde{W}_{{\alpha}}^T{\bf v}}{p}={c_\alpha}X^T\frac{(A_{c_\alpha}-\lambda_0I){\bf v}}{p}+\left[{c_\alpha}\beta,{A_{c_\alpha}\beta},\ldots,{A_{c_\alpha}^{n-1}\beta}\right]^T{\bf v}.
$$
It is straightforward to check that there exists an integer $s$ such that $\frac{p+s}{c_\alpha}\times {c_\alpha}\equiv1\pmod{p}$ and ${c_\alpha}\mid(p+s).$ Hence, $\frac{p+s}{c_\alpha}\equiv \frac{1}{c_\alpha}\pmod{p}.$ Notice that $\frac{\tilde{W}_{{\alpha}}^T{\bf v}}{p}$ is an integral vector. It follows that
\[\label{eq:3.01}
  \frac{\tilde{W}_{{\alpha}}^T{\bf v}}{p}\equiv X^T\frac{(A_{c_\alpha}-\lambda_0I){\bf v}}{p}+{\bf v}^T\beta\left[1,\frac{(p+s)\lambda_0}{c_\alpha},\ldots,\frac{(p+s)\lambda_0^{n-1}}{c_\alpha}\right]^T\pmod{p}.
\]
By Lemma \ref{lem3.2} and ${\bf v}^T{\bf v}\equiv0\pmod{p^2},$ one has ${\bf v}^T\frac{(A_{c_\alpha}-\lambda_0I){\bf v}}{p}\equiv0\pmod{p}.$ Together with ${\bf v}^T(A_{c_\alpha}-\lambda_0I)\equiv{\bf 0}\pmod{p}$ and $\rank_p{(A_{c_\alpha}-\lambda_0I)}=n-1,$ we obtain that there exists an integral vector ${\bf x}$ such that $\frac{(A_{c_\alpha}-\lambda_0I){\bf v}}{p}\equiv(A_{c_\alpha}-\lambda_0I){\bf x}\pmod{p}.$

In view of \eqref{eq:3.3}, one has
$$
  \tilde{W}_{{\alpha}}\equiv (A_{c_\alpha}-\lambda_0I)X+p\left[\beta,\frac{(p+s)A_{c_\alpha}\beta}{c_\alpha},\ldots,\frac{(p+s)A_{c_\alpha}^{n-1}\beta}{c_\alpha}\right]\equiv(A_{c_\alpha}-\lambda_0I)X\pmod{p}.
$$
Thus,
\[\label{eq:3.02}
  X^T\frac{(A_{c_\alpha}-\lambda_0I){\bf v}}{p}\equiv X^T(A_{c_\alpha}-\lambda_0I){\bf x}\equiv\tilde{W}_{{\alpha}}^T{\bf x}\pmod{p}.
\]
Recall that ${\bf v}\equiv(A_{c_\alpha}-\lambda_0I){\bf y}\pmod{p}.$ It is routine to check that
\begin{align*}
\frac{{\bf 1}^TA_{c_\alpha}{\bf y}}{c_\alpha}&\equiv\frac{(p+s)\lambda_0}{c_\alpha}{\bf 1}^T{\bf y}+\frac{p+s}{c_\alpha}{\bf 1}^T{\bf v}\equiv\frac{(p+s)\lambda_0}{c_\alpha}{\bf 1}^T{\bf y}\pmod{p},\\
\frac{{\bf 1}^TA_{c_\alpha}^2{\bf y}}{c_\alpha}&\equiv\frac{(p+s)\lambda_0}{c_\alpha}{\bf 1}^TA_{c_\alpha}{\bf y}+\frac{p+s}{c_\alpha}{\bf 1}^TA_{c_\alpha}{\bf v}\equiv\frac{(p+s)\lambda_0^2}{c_\alpha}{\bf 1}^T{\bf y}\pmod{p},
\end{align*}
\begin{align*}
    &\ \ \vdots\\
\frac{{\bf 1}^TA_{c_\alpha}^{n-1}{\bf y}}{c_\alpha}&\equiv\frac{(p+s)\lambda_0^{n-1}}{c_\alpha}{\bf 1}^T{\bf y}\pmod{p}.
\end{align*}
Hence,
\[\label{eq:3.03}
  \tilde{W}_{{\alpha}}^T{\bf y}\equiv {\bf 1}^T{\bf y}\left[1,\frac{(p+s)\lambda_0}{c_\alpha},\ldots,\frac{(p+s)\lambda_0^{n-1}}{c_\alpha}\right]^T \pmod{p}.
\]
Next, we show that ${\bf 1}^T{\bf y}\not\equiv0\pmod{p}.$ Otherwise, ${\bf 1}^T{\bf y}\equiv0\pmod{p}$, then $\tilde{W}_{{\alpha}}^T{\bf y}\equiv {\bf 0}\pmod{p}.$ Recall that $\tilde{W}_{{\alpha}}^T{\bf v}\equiv {\bf 0}\pmod{p}$ and $\rank_p(\tilde{W}_{{\alpha}})=n-1.$ It follows that ${\bf v}$ and ${\bf y}$ are linearly dependent over $\mathbb{F}_p.$ That is, there exist not all zero integers $m_1$ and $m_2$ such that $m_1{\bf v}+m_2{\bf y}={\bf 0}$ over $\mathbb{F}_p.$ Left multiplying both sides by $A_{c_\alpha}-\lambda_0I$ yields that $m_2{\bf v}\equiv{\bf 0}\pmod{p}.$ Notice that ${\bf v}\not\equiv{\bf 0}\pmod{p}.$ So, $m_2=0$ and therefore $m_1=0,$ a contradiction. Thus, ${\bf 1}^T{\bf y}\not\equiv0\pmod{p}.$ Then there exists an integer $t$ such that ${\bf v}^T\beta\equiv t{\bf 1}^T{\bf y}\pmod{p}.$

Together with \eqref{eq:3.01}-\eqref{eq:3.03}, we obtain that
$$
  \frac{\tilde{W}_{{\alpha}}^T{\bf v}}{p}\equiv \tilde{W}_{{\alpha}}^T{\bf x}+t\tilde{W}_{{\alpha}}^T{\bf y}\pmod{p},
$$
which is equivalent to
$$
  \tilde{W}_{{\alpha}}^T({\bf v}-p{\bf x}-tp{\bf y})\equiv{\bf 0}\pmod{p^2}.
$$
In view of Lemma \ref{lem2.4}, one has $p^2\mid \det\tilde{W}_{{\alpha}},$ a contradiction.

{\bf Case 2.}\ $\rank_p(A_{c_\alpha}-\lambda_0I)= n-2.$ In this case, we proceed by considering the following claim.
\begin{claim}\label{c1}
$\rank_p([A_{c_\alpha}-\lambda_0I,{\bf v}])=n-1.$
\end{claim}
\begin{proof}[\bf Proof of Claim \ref{c1}]\
Notice that $\rank_p([A_{c_\alpha}-\lambda_0I,{\bf v}])\geqslant \rank_p(A_{c_\alpha}-\lambda_0I)=n-2.$ Suppose that $\rank_p([A_{c_\alpha}-\lambda_0I,{\bf v}])=n-2.$ Then ${\bf v}$ can be written as a linear combination of the column vectors of $A_{c_\alpha}-\lambda_0I$ over $\mathbb{F}_p$. That is, there exist an integral vector ${\bf w}\,({\bf w}\not\equiv{\bf 0}\pmod{p})$ such that ${\bf v}\equiv (A_{c_\alpha}-\lambda_0I){\bf w}\pmod{p}.$ Recall that $A_{c_\alpha}{\bf v}\equiv \lambda_0{\bf v}\pmod{p}$ and $\bar{U}{\bf 1}=lU{\bf 1}=l{\bf 1}.$ Hence, for each positive integer $k,$
$$
  \frac{{\bf 1}^TA_{c_\alpha}^k{\bf w}}{c_\alpha}\equiv\frac{{\bf 1}^TA_{c_\alpha}^{k-1}{\bf v}}{c_\alpha}+\frac{\lambda_0{\bf 1}^TA_{c_\alpha}^{k-1}{\bf w}}{c_\alpha}\equiv\frac{\lambda_0^k{\bf 1}^T{\bf w}}{c_\alpha}\pmod{p}.
$$
Since $\rank_p(A_{c_\alpha}-\lambda_0I)= n-2,$ there exists an integral vector ${\bf y}$ which is linearly independent with ${\bf v}$ over $\mathbb{F}_p$, such that $(A_{c_\alpha}-\lambda_0I){\bf y}\equiv{\bf 0}\pmod{p}.$ It is routine to check that ${\bf 1}^T{\bf y}\not\equiv0\pmod{p}.$ In fact, if ${\bf 1}^T{\bf y}\equiv0\pmod{p},$ then $\tilde{W}_\alpha^T{\bf y}\equiv{\bf 0}\pmod{p},$ a contradiction to $\rank_p(\tilde{W}_\alpha)=n-1$.

Now, we show that ${\bf v},\,{\bf y}$ and ${\bf w}$ are linearly independent over $\mathbb{F}_p.$ Otherwise, there exist not all zeros integers $m_1,\,m_2$ and $m_3$ such that  $m_1{\bf v}+m_2{\bf y}+m_3{\bf w}={\bf 0}$ over $\mathbb{F}_p.$ Left multiplying both sides by $A_{c_\alpha}-\lambda_0I$ gives us $m_3{\bf w}\equiv{\bf 0}\pmod{p}.$ Notice that ${\bf w}\not\equiv{\bf 0}\pmod{p}.$ Hence, $m_3=0$ and so $m_1{\bf v}+m_2{\bf y}={\bf 0}$ over $\mathbb{F}_p,$ a contradiction. 

Let $\eta=({\bf 1}^T{\bf y}){\bf w}-({\bf 1}^T{\bf w}){\bf y}.$ Then $\eta\not\equiv{\bf 0}\pmod{p}$ and ${\bf 1}^T\eta\equiv0\pmod{p}.$ Moreover, for each integer  $k\geqslant 1,$
$$
  \frac{{\bf 1}^TA_{c_\alpha}^{k}\eta}{c_\alpha}\equiv({\bf 1}^T{\bf y})\frac{\lambda_0^k{\bf 1}^T{\bf w}}{c_\alpha}-({\bf 1}^T{\bf w})\frac{\lambda_0^k{\bf 1}^T{\bf y}}{c_\alpha}\equiv0\pmod{p}.
$$
It follows that $\tilde{W}_{{\alpha}}^T\eta\equiv{\bf 0}\pmod{p},$ a contradiction to the fact that $\rank_p(\tilde{W}_{{\alpha}})=n-1.$ This completes the proof of Claim \ref{c1}.
\end{proof}
Notice that ${\bf v}^T{\bf 1}\equiv0\pmod{p}$ and ${\bf v}^T[A_{c_\alpha}-\lambda_0I,{\bf v}]\equiv{\bf 0}\pmod{p}.$ By Claim \ref{c1}, there exist integral vectors ${\bf u},\,\beta$ and integer $f$ such that ${\bf 1}=(A_{c_\alpha}-\lambda_0I){\bf u}+f{\bf v}+p\beta.$ Therefore,
\[\label{eq:3.4}
  \tilde{W}_{{\alpha}}=(A_{c_\alpha}-\lambda_0I)X
    +f\left[{\bf v},\frac{A_{c_\alpha}{\bf v}}{c_\alpha},\ldots,\frac{A_{c_\alpha}^{n-1}{\bf v}}{c_\alpha}\right]+p\left[\beta,\frac{A_{c_\alpha}\beta}{c_\alpha},\ldots,\frac{A_{c_\alpha}^{n-1}\beta}{c_\alpha}\right],
\]
where $X=\left[{\bf u},\frac{A_{c_\alpha}{\bf u}}{c_\alpha},\ldots,\frac{A_{c_\alpha}^{n-1}{\bf u}}{c_\alpha}\right].$ Thus,
$$
  {c_\alpha}\frac{\tilde{W}_{{\alpha}}^T{\bf v}}{p}={c_\alpha}X^T\frac{(A_{c_\alpha}-\lambda_0I){\bf v}}{p}+\frac{f}{p}\left[{c_\alpha} {\bf v},A_{c_\alpha}{\bf v},\ldots,A_{c_\alpha}^{n-1}{\bf v}\right]^T{\bf v}+\left[{c_\alpha}\beta,{A_{c_\alpha}\beta},\ldots,{A_{c_\alpha}^{n-1}\beta}\right]^T{\bf v}.
$$
Notice that $\frac{\tilde{W}_{{\alpha}}^T{\bf v}}{p}$ is an integral vector. It follows that
$$
  \frac{\tilde{W}_{{\alpha}}^T{\bf v}}{p}\equiv X^T\frac{(A_{c_\alpha}-\lambda_0I){\bf v}}{p}+(\frac{f}{p}{\bf v}^T{\bf v}+{\bf v}^T\beta)\left[1,\frac{(p+s)\lambda_0}{c_\alpha},\ldots,\frac{(p+s)\lambda_0^{n-1}}{c_\alpha}\right]^T\pmod{p}.
$$
By Lemma \ref{lem3.2} and ${\bf v}^T{\bf v}\equiv0\pmod{p^2},$ one has ${\bf v}^T\frac{(A_{c_\alpha}-\lambda_0I){\bf v}}{p}\equiv0\pmod{p}.$ Together with Claim~\ref{c1} and ${\bf v}^T[A_{c_\alpha}-\lambda_0I,{\bf v}]\equiv{\bf 0}\pmod{p},$ there exists an integral vector ${\bf x}$ and an integer $m$ such that $\frac{(A_{c_\alpha}-\lambda_0I){\bf v}}{p}\equiv(A_{c_\alpha}-\lambda_0I){\bf x}+m{\bf v}\pmod{p}.$

Note that $(A_{c_\alpha}-\lambda_0 I){\bf y}\equiv{\bf 0}\pmod{p},$ here ${\bf y}$ is defined in the proof of Claim \ref{c1}. Thus,
$$
  \tilde{W}_{{\alpha}}^T{\bf y}=\left[{\bf 1},\frac{A_{c_\alpha}{\bf 1}}{c_\alpha},\ldots,\frac{A_{c_\alpha}^{n-1}{\bf 1}}{c_\alpha}\right]^T{\bf y}\equiv
  {\bf 1}^T{\bf y}\left[1,\frac{(p+s)\lambda_0}{c_\alpha},\ldots,\frac{(p+s)\lambda_0^{n-1}}{c_\alpha}\right]^T\pmod{p}.
$$
Recall that ${\bf 1}^T{\bf y}\not\equiv0\pmod{p}.$ Hence, there exists an integer $t$ satisfying $\frac{f}{p}{\bf v}^T{\bf v}+{\bf v}^T\beta+m{\bf u}^T{\bf v}-f{\bf v}^T{\bf x}\equiv t{\bf 1}^T{\bf y}\pmod{p}.$ It follows from \eqref{eq:3.4} that $\tilde{W}_{{\alpha}}^T\equiv X^T{(A_{c_\alpha}-\lambda_0I)}+f\left[{\bf v},\frac{A_{c_\alpha}{\bf v}}{c_\alpha},\ldots,\frac{A_{c_\alpha}^{n-1}{\bf v}}{c_\alpha}\right]\pmod{p}.$ Hence
\begin{align*}
    \frac{\tilde{W}_{{\alpha}}^T{\bf v}}{p}&\equiv X^T\frac{(A_{c_\alpha}-\lambda_0I){\bf v}}{p}+(\frac{f}{p}{\bf v}^T{\bf v}+{\bf v}^T\beta)\left[1,\frac{(p+s)\lambda_0}{c_\alpha},\ldots,\frac{(p+s)\lambda_0^{n-1}}{c_\alpha}\right]^T\\
    &\equiv X^T(A_{c_\alpha}-\lambda_0I){\bf x}+mX^T{\bf v}+(\frac{f}{p}{\bf v}^T{\bf v}+{\bf v}^T\beta)\left[1,\frac{(p+s)\lambda_0}{c_\alpha},\ldots,\frac{(p+s)\lambda_0^{n-1}}{c_\alpha}\right]^T\\
    &\equiv\tilde{W}_{{\alpha}}^T{\bf x}+(\frac{f}{p}{\bf v}^T{\bf v}+{\bf v}^T\beta+m{\bf u}^T{\bf v}-f{\bf v}^T{\bf x})\left[1,\frac{(p+s)\lambda_0}{c_\alpha},\ldots,\frac{(p+s)\lambda_0^{n-1}}{c_\alpha}\right]^T\\
    &\equiv\tilde{W}_{{\alpha}}^T{\bf x}+t{\bf 1}^T{\bf y}\left[1,\frac{(p+s)\lambda_0}{c_\alpha},\ldots,\frac{(p+s)\lambda_0^{n-1}}{c_\alpha}\right]^T\\
    &\equiv\tilde{W}_{{\alpha}}^T{\bf x}+t\tilde{W}_{{\alpha}}^T{\bf y}\pmod{p}.
\end{align*}
Therefore,
$$
  \tilde{W}_{{\alpha}}^T({\bf v}-p{\bf x}-tp{\bf y})\equiv{\bf 0}\pmod{p^2}.
$$
In view of Lemma \ref{lem2.4}, one has $p^2\mid \det\tilde{W}_{{\alpha}},$ a contradiction.

{\bf Case 3.}\ $\rank_p(A_{c_\alpha}-\lambda_0I)\leqslant n-3.$ In this case, there exists three linearly independent vectors, say ${\bf v},\,{\bf w}$ and ${\bf y},$ such that $(A_{c_\alpha}-\lambda_0I){\bf v}=(A_{c_\alpha}-\lambda_0I){\bf w}=(A_{c_\alpha}-\lambda_0I){\bf y}={\bf 0}$ over $\mathbb{F}_p.$ It is straightforward to check that ${\bf 1}^T{\bf w}\not\equiv0\pmod{p}$ and ${\bf 1}^T{\bf y}\not\equiv0\pmod{p}.$ Let $\zeta=({\bf 1}^T{\bf y}){\bf w}-({\bf 1}^T{\bf w}){\bf y}.$ Then $\zeta\not\equiv{\bf 0}\pmod{p}$ and ${\bf 1}^T\zeta\equiv0\pmod{p}.$ Moreover, for each integer $k\geqslant 1,$
$$
  \frac{{\bf 1}^TA_{c_\alpha}^{k}\zeta}{c_\alpha}\equiv({\bf 1}^T{\bf y})\frac{\lambda_0^k{\bf 1}^T{\bf w}}{c_\alpha}-({\bf 1}^T{\bf w})\frac{\lambda_0^k{\bf 1}^T{\bf y}}{c_\alpha}\equiv0\pmod{p}.
$$
It follows that $\tilde{W}_{{\alpha}}^T\zeta\equiv{\bf 0}\pmod{p},$ a contradiction to the fact that $\rank_p(\tilde{W}_{{\alpha}})=n-1.$

Combining Cases 1-3, Theorem \ref{thm3.1} follows immediately.
\end{proof}
\subsection{\normalsize Proof of Theorem \ref{thm3.3}}
In this subsection, we present the proof of Theorem \ref{thm3.3}. Before doing so, we need the following lemmas.
\begin{lem}\label{lem3.4}
Let $G$ be a graph and let $\alpha$ be in $[0,1)$. Then $\rank_2(\tilde{W}_{{\alpha}})\leqslant \lceil\frac{n}{2}\rceil.$
\end{lem}
\begin{proof}
We proceed by distinguishing the parity on $n$.

{\bf Case 1.}\ $n$ is even. In view of Lemma \ref{lem2.5}, one obtains that
\begin{equation}\label{eq:4}
    \tilde{W}_{{\alpha}}^T\tilde{W}_{{\alpha}}=\left(
                                               \begin{array}{cccc}
                                                 {\bf 1}^T{\bf 1} & \frac{{\bf 1}^TA_{c_\alpha}{\bf 1}}{c_\alpha} & \cdots & \frac{{\bf 1}^TA_{c_\alpha}^{n-1}{\bf 1}}{c_\alpha} \\
                                                 \frac{{\bf 1}^TA_{c_\alpha}{\bf 1}}{c_\alpha} & \frac{{\bf 1}^TA_{c_\alpha}^2{\bf 1}}{{c_\alpha^2}} & \cdots & \frac{{\bf 1}^TA_{c_\alpha}^{n}{\bf 1}}{{c_\alpha^2}} \\
                                                 \vdots & \vdots & \ddots & \vdots \\
                                                 \frac{{\bf 1}^TA_{c_\alpha}^{n-1}{\bf 1}}{c_\alpha} & \frac{{\bf 1}^TA_{c_\alpha}^{n}{\bf 1}}{{c_\alpha^2}} & \cdots & \frac{{\bf 1}^TA_{c_\alpha}^{2n-2}{\bf 1}}{{c_\alpha^2}} \\
                                               \end{array}
                                             \right)\equiv{\bf 0}\pmod{2}.
\end{equation}
Therefore, $2\,\rank_2(\tilde{W}_{{\alpha}})=\rank_2(\tilde{W}_{{\alpha}}^T)+\rank_2(\tilde{W}_{{\alpha}})\leqslant n$. It follows that $\rank_2(\tilde{W}_{{\alpha}})\leqslant \frac{n}{2}=\lceil\frac{n}{2}\rceil.$

{\bf Case 2.}\ $n$ is odd. Let $\bar{W}_{{\alpha}}:=\bar{W}_{{\alpha}}(G)=[2\times{\bf 1},\frac{A_{c_\alpha}{\bf 1}}{c_\alpha},\ldots,\frac{A_{c_\alpha}^{n-1}{\bf 1}}{c_\alpha}].$ Based on Lemma \ref{lem2.5}, we obtain that
\begin{equation}\label{eq:5}
    \tilde{W}_{{\alpha}}^T\bar{W}_{{\alpha}}=\left(
                                               \begin{array}{cccc}
                                                 2\times{\bf 1}^T{\bf 1} & \frac{{\bf 1}^TA_{c_\alpha}{\bf 1}}{c_\alpha} & \cdots & \frac{{\bf 1}^TA_{c_\alpha}^{n-1}{\bf 1}}{c_\alpha} \\
                                                 \frac{2\times{\bf 1}^TA_{c_\alpha}{\bf 1}}{c_\alpha} & \frac{{\bf 1}^TA_{c_\alpha}^2{\bf 1}}{{c_\alpha^2}} & \cdots & \frac{{\bf 1}^TA_{c_\alpha}^{n}{\bf 1}}{{c_\alpha^2}} \\
                                                 \vdots & \vdots & \ddots & \vdots \\
                                                 \frac{2\times{\bf 1}^TA_{c_\alpha}^{n-1}{\bf 1}}{c_\alpha} & \frac{{\bf 1}^TA_{c_\alpha}^{n}{\bf 1}}{{c_\alpha^2}} & \cdots & \frac{{\bf 1}^TA_{c_\alpha}^{2n-2}{\bf 1}}{{c_\alpha^2}} \\
                                               \end{array}
                                             \right)\equiv{\bf 0}\pmod{2}.
\end{equation}
It is easy to see that $\rank_2(\bar{W}_{{\alpha}})\geqslant \rank_2(\tilde{W}_{{\alpha}})-1$ and $\rank_2(\tilde{W}_{{\alpha}}^T)+\rank_2(\bar{W}_{{\alpha}})\leqslant n.$ Hence, $\rank_2(\tilde{W}_{{\alpha}})\leqslant \frac{n+1}{2}=\lceil\frac{n}{2}\rceil.$

Combining Cases 1-2, we complete the proof.
\end{proof}
\begin{lem}\label{lem3.5}
Let $G\in \mathfrak{F}_n$ and $\alpha\in[0,1).$ Then the SNF of $\tilde{W}_{{\alpha}}$ is
$$
S={\rm diag}(\underbrace{1,\ldots,1}_{\lceil\frac{n}{2}\rceil},\underbrace{2,\ldots,2,2B}_{\lfloor\frac{n}{2}\rfloor}),
$$
where $B$ is an odd and square-free integer. Furthermore, $\rank_2(\tilde{W}_{{\alpha}})=\lceil\frac{n}{2}\rceil.$
\end{lem}
\begin{proof}
Notice that $\frac{\det \tilde{W}_{{\alpha}}}{2^{\lfloor\frac{n}{2}\rfloor}}$ is odd and square-free. Then $\det \tilde{W}_{{\alpha}}=\pm 2^{\lfloor\frac{n}{2}\rfloor}p_1p_2\ldots p_s,$ here $p_i$ denotes an odd prime and $p_i\neq p_j$ for $1\leqslant i<j\leqslant s.$ Therefore, the SNF of $\tilde{W}_{{\alpha}}$ is
$$
S={\rm diag}(1,\ldots,1,2^{l_1},\ldots,2^{l_{t-1}},2^{l_t}B),
$$
where $B=p_1p_2\ldots p_s.$ In view of Lemma \ref{lem3.4}, one has $\rank_2(\tilde{W}_{{\alpha}})\leqslant \lceil\frac{n}{2}\rceil$, which is equivalent to  $n-t\leqslant \lceil\frac{n}{2}\rceil.$ It follows that $t\geqslant \lfloor\frac{n}{2}\rfloor.$ Note that $\det(\tilde{W}_{{\alpha}})=\pm \det(S).$ Hence, $l_1+l_2+\cdots+l_t=\lfloor\frac{n}{2}\rfloor.$ Thus, $l_1=l_2=\cdots=l_t=1$ and $t=\lfloor\frac{n}{2}\rfloor.$ So, $\rank_2(\tilde{W}_{{\alpha}})=\lceil\frac{n}{2}\rceil.$ This completes the proof.
\end{proof}

For the ease of presentation, we need the following notations. For a graph $G$ on $n$ vertices, let $\hat{W}_{{\alpha}}(G)$ be the matrix defined as follows:
\begin{equation*}
     \hat{W}_{{\alpha}}:=\hat{W}_{{\alpha}}(G)=
    \left\{
    \begin{aligned}
        &\left[{\bf 1},\frac{A_{c_\alpha}{\bf 1}}{c_\alpha},\ldots,\frac{A_{c_\alpha}^{\frac{n}{2}-1}{\bf 1}}{c_\alpha}\right],&\ \ \textrm{if $n$ is even;}\\
        &\left[\frac{A_{c_\alpha}{\bf 1}}{c_\alpha},\frac{A_{c_\alpha}^2{\bf 1}}{c_\alpha},\ldots,\frac{A_{c_\alpha}^{\frac{n-1}{2}}{\bf 1}}{c_\alpha}\right],&\ \ \textrm{if $n$ is odd.}
    \end{aligned}
    \right.
\end{equation*}
\begin{lem}\label{lem3.7}
Let $G\in \mathfrak{F}_n$. Then $\rank_2(\hat{W}_{{\alpha}})=\lfloor\frac{n}{2}\rfloor$ for $\alpha\in[0,1).$
\end{lem}
\begin{proof}
By Lemma \ref{lem3.5}, one has $\rank_2(\tilde{W}_{{\alpha}})=\lceil\frac{n}{2}\rceil.$ Let $t:=\lceil\frac{n}{2}\rceil.$ In order to prove our result, it suffices to show that the first $t$ columns of $\tilde{W}_{{\alpha}}$ are linearly independent over $\mathbb{F}_2.$ Suppose, to the contrary, that ${\bf 1},\frac{A_{c_\alpha}{\bf 1}}{c_\alpha},\ldots,\frac{A_{c_\alpha}^{t-1}{\bf 1}}{c_\alpha}$ are linearly dependent over $\mathbb{F}_2.$ That is to say, there exist not all zero integers $b_0,b_1,\ldots,b_{t-1}\in \mathbb{F}_2$ such that
$$
 b_0{\bf 1}+b_1\frac{A_{c_\alpha}{\bf 1}}{c_\alpha}+\cdots+b_{t-1}\frac{A_{c_\alpha}^{t-1}{\bf 1}}{c_\alpha}\equiv{\bf 0}\pmod{2}.
$$
Let $m=\max\{i:0\leqslant i\leqslant t-1,b_i\neq0\}.$ Clearly,  $0<m\leqslant t-1.$ Hence
$$
  \frac{A_{c_\alpha}^{m}{\bf 1}}{c_\alpha}=-b_m^{-1}b_0{\bf 1}-b_m^{-1}b_1\frac{A_{c_\alpha}{\bf 1}}{c_\alpha}-\cdots-b_m^{-1}b_{m-1}\frac{A_{c_\alpha}^{m-1}{\bf 1}}{c_\alpha}\ \text{over}\ \mathbb{F}_2.
$$
It follows that $\frac{A_{c_\alpha}^{m}{\bf 1}}{c_\alpha}$ can be written as a linear combination of ${\bf 1},\frac{A_{c_\alpha}{\bf 1}}{c_\alpha},\ldots, \frac{A_{c_\alpha}^{m-1}{\bf 1}}{c_\alpha}$ over $\mathbb{F}_2.$ Then
$$
  \frac{A_{c_\alpha}^{m}{\bf 1}}{c_\alpha}=-b_m^{-1}b_0{\bf 1}-b_m^{-1}b_1\frac{A_{c_\alpha}{\bf 1}}{c_\alpha}-\cdots-b_m^{-1}b_{m-1}\frac{A_{c_\alpha}^{m-1}{\bf 1}}{c_\alpha}+2\beta\ \text{over}\ \mathbb{Z}
$$
for some integral vector $\beta.$ Moreover,
$$
  \frac{A_{c_\alpha}^{m+1}{\bf 1}}{c_\alpha}=-{c_\alpha}b_m^{-1}b_0\frac{A_{c_\alpha}{\bf 1}}{c_\alpha}-b_m^{-1}b_1\frac{A_{c_\alpha}^2{\bf 1}}{c_\alpha}-\cdots-b_m^{-1}b_{m-1}\frac{A_{c_\alpha}^{m}{\bf 1}}{c_\alpha}+2A_{c_\alpha}\beta\ \text{over}\ \mathbb{Z}.
$$
Hence
$$
  \frac{A_{c_\alpha}^{m+1}{\bf 1}}{c_\alpha}{\color{blue}=}-{c_\alpha}b_m^{-1}b_0\frac{A_{c_\alpha}{\bf 1}}{c_\alpha}-b_m^{-1}b_1\frac{A_{c_\alpha}^2{\bf 1}}{c_\alpha}-\cdots-b_m^{-1}b_{m-1}\frac{A_{c_\alpha}^{m}{\bf 1}}{c_\alpha}\ \text{over}\ \mathbb{F}_2,
$$
i.e., $\frac{A_{c_\alpha}^{m+1}{\bf 1}}{c_\alpha}$ can be written as a linear combination of ${\bf 1},\frac{A_{c_\alpha}{\bf 1}}{c_\alpha},\ldots, \frac{A_{c_\alpha}^{m-1}{\bf 1}}{c_\alpha}$ over $\mathbb{F}_2.$

By a similar discussion, we can show that for each integer $i\geqslant 0,$ the vector $\frac{A_{c_\alpha}^{m+i}{\bf 1}}{c_\alpha}$ can be written as a linear combination of ${\bf 1},\frac{A_{c_\alpha}{\bf 1}}{c_\alpha},\ldots, \frac{A_{c_\alpha}^{m-1}{\bf 1}}{c_\alpha}$ over $\mathbb{F}_2.$
Therefore, $\rank_2(\tilde{W}_{{\alpha}})\leqslant m\leqslant t-1,$ a contradiction.
This completes the proof.
\end{proof}
For a graph $G,$ let $\tilde{W}'_{{\alpha}}:=\tilde{W}'_{{\alpha}}(G)=\left[{\bf 1},\frac{A_{c_\alpha}^2{\bf 1}}{c_\alpha},\ldots,\frac{A_{c_\alpha}^{2n-2}{\bf 1}}{c_\alpha}\right],$ and let
\begin{equation*}
    \hat{W}'_{{\alpha}}:=\hat{W}'_{{\alpha}}(G)=
    \left\{
    \begin{aligned}
        &\left[{\bf 1},\frac{A_{c_\alpha}^2{\bf 1}}{c_\alpha},\ldots,\frac{A_{c_\alpha}^{n-2}{\bf 1}}{c_\alpha}\right],&\ \ \textrm{if $n$ is even;}\\
        &\left[\frac{A_{c_\alpha}^2{\bf 1}}{c_\alpha},\frac{A_{c_\alpha}^4{\bf 1}}{c_\alpha},\ldots,\frac{A_{c_\alpha}^{n-1}{\bf 1}}{c_\alpha}\right],&\ \ \textrm{if $n$ is odd.}
    \end{aligned}
    \right.
\end{equation*}
\begin{lem}\label{lem3.8}
Let $G\in \mathfrak{F}_n$. Then $\rank_2\left(\frac{\tilde{W}_{{\alpha}}^T\hat{W}'_{{\alpha}}}{2}\right)=\lfloor\frac{n}{2}\rfloor$ for $\alpha\in[0,1).$
\end{lem}
\begin{proof}
For even $n,$ in view of Lemma \ref{lem3.5}, one has $\det \left(\frac{\tilde{W}_{{\alpha}}^T\tilde{W}_{{\alpha}}}{2}\right)=\frac{(2^{\lfloor\frac{n}{2}\rfloor}B)^2}{2^n}=B^2,$ where $B$ is defined in Lemma \ref{lem3.5}. Notice that $B$ is odd. Hence, $\rank_2 \left(\frac{\tilde{W}_{{\alpha}}^T\tilde{W}_{{\alpha}}}{2}\right)=n.$ Therefore, $\rank_2\left(\frac{\tilde{W}_{{\alpha}}^T\hat{W}'_{{\alpha}}}{2}\right)$ equals the number of columns of $\hat{W}'_{{\alpha}},$ as desired. 

For odd $n,$ let $\bar{W}_{{\alpha}}:=\bar{W}_{{\alpha}}(G)=[2\times {\bf 1},\frac{A_{c_\alpha}{\bf 1}}{c_\alpha},\ldots,\frac{A_{c_\alpha}^{n-1}{\bf 1}}{c_\alpha}]$ be the matrix defined as Lemma~\ref{lem3.4}. In view of Lemma \ref{lem3.5}, one has $\det \left(\frac{\tilde{W}_{{\alpha}}^T\bar{W}_{{\alpha}}}{2}\right)=B^2.$  It follows that $\rank_2 \left(\frac{\tilde{W}_{{\alpha}}^T\bar{W}_{{\alpha}}}{2}\right)=n.$ Hence,  $\rank_2\left(\frac{\tilde{W}_{{\alpha}}^T\hat{W}'_{{\alpha}}}{2}\right)$ is equal to the number of columns of $\hat{W}'_{{\alpha}},$ as desired. This completes the proof. 
\end{proof}
Now, we are ready to present the proof Theorem \ref{thm3.3}.
\begin{proof}[\bf Proof of Theorem \ref{thm3.3}]\
Suppose on the contrary that $l$ is even. By Lemma \ref{lem3.2}, there exists a column ${\bf v}\,({\bf v}\not\equiv{\bf 0}\pmod{2})$ of $lU$ satisfying that $\tilde{W}_{{\alpha}}^T{\bf v}\equiv{\bf 0}\pmod{2}$ and ${\bf v}^TA_{c_\alpha}^k{\bf v}\equiv0\pmod{4}$ for each integer $k\geqslant 0.$ In view of \eqref{eq:4}, \eqref{eq:5} and Lemma \ref{lem3.7}, we obtain that ${\bf v}$ can be written as a linear combination of the column vectors of $\hat{W}_{{\alpha}}$ over $\mathbb{F}_2,$ that is, there exist integral vectors ${\bf u}\,({\bf u}\not\equiv{\bf 0}\pmod{2})$ and $\beta$ such that ${\bf v}=\hat{W}_{{\alpha}}{\bf u}+2\beta.$ Therefore, for each integer $k\geqslant 0,$

\begin{align*}
    {\bf v}^TA_{c_\alpha}^k{\bf v}=(\hat{W}_{{\alpha}}{\bf u}+2\beta)^TA_{c_\alpha}^k(\hat{W}_{{\alpha}}{\bf u}+2\beta)\equiv {\bf u}^T\hat{W}_{{\alpha}}^TA_{c_\alpha}^k\hat{W}_{{\alpha}}{\bf u}\equiv0\pmod{4}.
\end{align*}
We proceed by distinguishing the parity of $n$.

{\bf Case 1.}\ $n$ is even.  In view of Lemma \ref{lem2.5}, one has for each integer $k\geqslant 0,$
\begin{equation*}
    \hat{W}_{{\alpha}}^TA_{c_\alpha}^k\hat{W}_{{\alpha}}=\left[
                                                    \begin{array}{cccc}
                                                      {\bf 1}^TA_{c_\alpha}^k{\bf 1} & \frac{{\bf 1}^TA_{c_\alpha}^{1+k}{\bf 1}}{c_\alpha} & \cdots & \frac{{\bf 1}^TA_{c_\alpha}^{\frac{n}{2}-1+k}{\bf 1}}{c_\alpha} \\
                                                      \frac{{\bf 1}^TA_{c_\alpha}^{1+k}{\bf 1}}{c_\alpha} & \frac{{\bf 1}^TA_{c_\alpha}^{2+k}{\bf 1}}{{c_\alpha^2}} & \cdots & \frac{{\bf 1}^TA_{c_\alpha}^{\frac{n}{2}+k}{\bf 1}}{{c_\alpha^2}} \\
                                                      \vdots & \vdots & \ddots & \vdots \\
                                                      \frac{{\bf 1}^TA_{c_\alpha}^{\frac{n}{2}-1+k}{\bf 1}}{c_\alpha} & \frac{{\bf 1}^TA_{c_\alpha}^{\frac{n}{2}+k}{\bf 1}}{{c_\alpha^2}} & \cdots & \frac{{\bf 1}^TA_{c_\alpha}^{n-2+k}{\bf 1}}{{c_\alpha^2}} \\
                                                    \end{array}
                                                  \right]\equiv{\bf 0}\pmod{2}.
\end{equation*}
Put $M:=(M_{ij})_{\frac{n}{2}\times\frac{n}{2}}=\hat{W}_{{\alpha}}^TA_{c_\alpha}^k\hat{W}_{{\alpha}}$ and ${\bf u}=[u_1,u_2,\ldots,u_{\frac{n}{2}}]^T.$ 
Hence, for each integer $k\geqslant 0,$
\begin{align}
    {\bf u}^T\hat{W}_{{\alpha}}^TA_{c_\alpha}^k\hat{W}_{{\alpha}}{\bf u}&=\sum_{1\leqslant i\leqslant j\leqslant {\frac{n}{2}}}M_{ij}u_iu_j=\sum_{1\leqslant i\leqslant {\frac{n}{2}}}M_{ii}u_i^2+2\sum_{1\leqslant i<j\leqslant {\frac{n}{2}}}M_{ij}u_iu_j\notag\\
    &\equiv ({\bf 1}^TA_{c_\alpha}^k{\bf 1})u_1+\frac{{\bf 1}^TA_{c_\alpha}^{2+k}{\bf 1}}{{c_\alpha^2}}u_2+\cdots+\frac{{\bf 1}^TA_{c_\alpha}^{n-2+k}{\bf 1}}{{c_\alpha^2}}u_{\frac{n}{2}}\label{eq:3.09}\\
    &=\left[{\bf 1}^TA_{c_\alpha}^k{\bf 1},\frac{{\bf 1}^TA_{c_\alpha}^{2+k}{\bf 1}}{{c_\alpha^2}},\ldots,\frac{{\bf 1}^TA_{c_\alpha}^{n-2+k}{\bf 1}}{{c_\alpha^2}}\right]{\bf u}\notag\\
    &\equiv0\pmod{4},\notag
\end{align}
the congruence equation in \eqref{eq:3.09} follows from the fact that $u_i^2\equiv u_i\pmod{2}$ for any $1\leqslant i\leqslant \frac{n}{2}.$ Therefore,
\[\label{eq:3.5}
\left[\frac{{\bf 1}^TA_{c_\alpha}^k{\bf 1}}{2},\frac{{\bf 1}^TA_{c_\alpha}^{2+k}{\bf 1}}{2{c_\alpha^2}},\ldots,\frac{{\bf 1}^TA_{c_\alpha}^{n-2+k}{\bf 1}}{2{c_\alpha^2}}\right]{\bf u}\equiv0\pmod{2}.
\]

If $c_\alpha$ is even, then we define
\begin{align*}
    M_1:&=\left[
          \begin{array}{ccccc}
            \frac{{\bf 1}^T{\bf 1}}{2} & \frac{{\bf 1}^TA_{c_\alpha}^2{\bf 1}}{2{c_\alpha^2}} & \frac{{\bf 1}^TA_{c_\alpha}^4{\bf 1}}{2{c_\alpha^2}} & \cdots & \frac{{\bf 1}^TA_{c_\alpha}^{n-2}{\bf 1}}{2{c_\alpha^2}} \\
            0 & \frac{{\bf 1}^TA_{c_\alpha}^{3}{\bf 1}}{2{c_\alpha^2}} & \frac{{\bf 1}^TA_{c_\alpha}^{5}{\bf 1}}{2{c_\alpha^2}} & \cdots & \frac{{\bf 1}^TA_{c_\alpha}^{n-1}{\bf 1}}{2{c_\alpha^2}} \\
            0 & \frac{{\bf 1}^TA_{c_\alpha}^{4}{\bf 1}}{2{c_\alpha^2}} & \frac{{\bf 1}^TA_{c_\alpha}^{6}{\bf 1}}{2{c_\alpha^2}} & \cdots & \frac{{\bf 1}^TA_{c_\alpha}^{n}{\bf 1}}{2{c_\alpha^2}} \\
            \vdots & \vdots & \vdots & \ddots & \vdots \\
            0 & \frac{{\bf 1}^TA_{c_\alpha}^{n+1}{\bf 1}}{2{c_\alpha^2}} & \frac{{\bf 1}^TA_{c_\alpha}^{n+3}{\bf 1}}{2{c_\alpha^2}} & \cdots & \frac{{\bf 1}^TA_{c_\alpha}^{2n-3}{\bf 1}}{2{c_\alpha^2}} \\
          \end{array}
        \right]\\
        &\equiv\left[
          \begin{array}{ccccc}
            \frac{{\bf 1}^T{\bf 1}}{2} & \frac{{\bf 1}^TA_{c_\alpha}^2{\bf 1}}{2{c_\alpha^2}} & \frac{{\bf 1}^TA_{c_\alpha}^4{\bf 1}}{2{c_\alpha^2}} & \cdots & \frac{{\bf 1}^TA_{c_\alpha}^{n-2}{\bf 1}}{2{c_\alpha^2}} \\
            \frac{{\bf 1}^TA_{c_\alpha}{\bf 1}}{2} & \frac{{\bf 1}^TA_{c_\alpha}^{3}{\bf 1}}{2{c_\alpha^2}} & \frac{{\bf 1}^TA_{c_\alpha}^{5}{\bf 1}}{2{c_\alpha^2}} & \cdots & \frac{{\bf 1}^TA_{c_\alpha}^{n-1}{\bf 1}}{2{c_\alpha^2}} \\
            \frac{{\bf 1}^TA_{c_\alpha}^2{\bf 1}}{2} & \frac{{\bf 1}^TA_{c_\alpha}^{4}{\bf 1}}{2{c_\alpha^2}} & \frac{{\bf 1}^TA_{c_\alpha}^{6}{\bf 1}}{2{c_\alpha^2}} & \cdots & \frac{{\bf 1}^TA_{c_\alpha}^{n}{\bf 1}}{2{c_\alpha^2}} \\
            \vdots & \vdots & \vdots & \ddots & \vdots \\
            \frac{{\bf 1}^TA_{c_\alpha}^{n-1}{\bf 1}}{2} & \frac{{\bf 1}^TA_{c_\alpha}^{n+1}{\bf 1}}{2{c_\alpha^2}} & \frac{{\bf 1}^TA_{c_\alpha}^{n+3}{\bf 1}}{2{c_\alpha^2}} & \cdots & \frac{{\bf 1}^TA_{c_\alpha}^{2n-3}{\bf 1}}{2{c_\alpha^2}} \\
          \end{array}
        \right]\pmod{2},
\end{align*}
the second congruence equation follows from Lemma \ref{lem2.5}. In view of \eqref{eq:3.5}, one has $M_1{\bf u}\equiv{\bf 0}\pmod{2}.$ Furthermore, we define
\begin{align*}
    M_2:&=\left[
          \begin{array}{ccccc}
            \frac{{\bf 1}^T{\bf 1}}{2} & 0 & 0 & \cdots & 0 \\
            \frac{{\bf 1}^TA_{c_\alpha}{\bf 1}}{2{c_\alpha}} & \frac{{\bf 1}^TA_{c_\alpha}^{3}{\bf 1}}{2{c_\alpha^2}} & \frac{{\bf 1}^TA_{c_\alpha}^{5}{\bf 1}}{2{c_\alpha^2}} & \cdots & \frac{{\bf 1}^TA_{c_\alpha}^{n-1}{\bf 1}}{2{c_\alpha^2}} \\
            0 & \frac{{\bf 1}^TA_{c_\alpha}^{4}{\bf 1}}{2{c_\alpha^2}} & \frac{{\bf 1}^TA_{c_\alpha}^{6}{\bf 1}}{2{c_\alpha^2}} & \cdots & \frac{{\bf 1}^TA_{c_\alpha}^{n}{\bf 1}}{2{c_\alpha^2}} \\
            \vdots & \vdots & \vdots & \ddots & \vdots \\
            0 & \frac{{\bf 1}^TA_{c_\alpha}^{n+1}{\bf 1}}{2{c_\alpha^2}} & \frac{{\bf 1}^TA_{c_\alpha}^{n+3}{\bf 1}}{2{c_\alpha^2}} & \cdots & \frac{{\bf 1}^TA_{c_\alpha}^{2n-3}{\bf 1}}{2{c_\alpha^2}} \\
          \end{array}
        \right]\\
        &\equiv\left[
          \begin{array}{ccccc}
            \frac{{\bf 1}^T{\bf 1}}{2} & \frac{{\bf 1}^TA_{c_\alpha}^2{\bf 1}}{2{c_\alpha}} & \frac{{\bf 1}^TA_{c_\alpha}^4{\bf 1}}{2{c_\alpha}} & \cdots & \frac{{\bf 1}^TA_{c_\alpha}^{n-2}{\bf 1}}{2{c_\alpha}} \\
            \frac{{\bf 1}^TA_{c_\alpha}{\bf 1}}{2{c_\alpha}} & \frac{{\bf 1}^TA_{c_\alpha}^{3}{\bf 1}}{2{c_\alpha^2}} & \frac{{\bf 1}^TA_{c_\alpha}^{5}{\bf 1}}{2{c_\alpha^2}} & \cdots & \frac{{\bf 1}^TA_{c_\alpha}^{n-1}{\bf 1}}{2{c_\alpha^2}} \\
            \frac{{\bf 1}^TA_{c_\alpha}^2{\bf 1}}{2{c_\alpha}} & \frac{{\bf 1}^TA_{c_\alpha}^{4}{\bf 1}}{2{c_\alpha^2}} & \frac{{\bf 1}^TA_{c_\alpha}^{6}{\bf 1}}{2{c_\alpha^2}} & \cdots & \frac{{\bf 1}^TA_{c_\alpha}^{n}{\bf 1}}{2{c_\alpha^2}} \\
            \vdots & \vdots & \vdots & \ddots & \vdots \\
            \frac{{\bf 1}^TA_{c_\alpha}^{n-1}{\bf 1}}{2{c_\alpha}} & \frac{{\bf 1}^TA_{c_\alpha}^{n+1}{\bf 1}}{2{c_\alpha^2}} & \frac{{\bf 1}^TA_{c_\alpha}^{n+3}{\bf 1}}{2{c_\alpha^2}} & \cdots & \frac{{\bf 1}^TA_{c_\alpha}^{2n-3}{\bf 1}}{2{c_\alpha^2}} \\
          \end{array}
        \right]
        =\frac{\tilde{W}_{{\alpha}}^T\hat{W}'_{{\alpha}}}{2}\pmod{2},
\end{align*}
the second congruence equation follows from Lemma \ref{lem2.5}.

By Lemma \ref{lem3.7}, we obtain that $\rank_2\left(\frac{\tilde{W}_{{\alpha}}^T\hat{W}'_{{\alpha}}}{2}\right)=\frac{n}{2},$ that is, the column rank of $M_2$ is full over $\mathbb{F}_2.$ Hence, the last $\frac{n}{2}-1$ columns of $M_2$ are linearly independent over $\mathbb{F}_2.$ Therefore, the last $\frac{n}{2}-1$ columns of $M_1$ are linearly independent over $\mathbb{F}_2.$ It follows that $\rank_2(M_1)\geqslant \frac{n}{2}-1.$ 

By applying \eqref{eq:3.5} with $k=0$, we have $\frac{{\bf 1}^T{\bf 1}}{2}\equiv0\pmod{2}.$ Therefore, $\rank_2(M_1)= \frac{n}{2}-1.$ Then the solution space of $M_1{\bf u}\equiv{\bf 0}\pmod{2}$ has dimension $1$. It is routine to check that it is spanned by ${\bf u}\equiv[1,0,\ldots,0]^T\pmod{2}.$ Hence, ${\bf v}\equiv\hat{W}_{{\alpha}}{\bf u}\equiv {\bf 1}\pmod{2}.$ In view of Lemmas \ref{lem3.6} and \ref{lem3.5}, one has $l\mid2B.$ It follows from Theorem~\ref{thm3.1} that $l\mid2$ and therefore $l=2.$ Recall that $U$ is orthogonal and $U{\bf 1}={\bf 1}.$ Based on the choice of ${\bf v}$ and $\frac{{\bf 1}^T{\bf 1}}{2}\equiv0\pmod{2},$ we know that $n=4.$  However, there are only six connected graphs on $4$ vertices. It is routine to check that each of them is DG$A_\alpha$S. Hence $l=1,$ a contradiction.

If $c_\alpha=1,$ then by \eqref{eq:3.5}, one has $M_2{\bf u}\equiv{\bf 0}\pmod{2}.$ Together with Lemma \ref{lem3.7}, we obtain ${\bf u}\equiv{\bf 0}\pmod{2},$ a contradiction.

{\bf Case 2.}\ $n$ is odd. By a similar discussion as the proof for even $n$, we can show that for each integer $k\geqslant 0,$
\[\label{eq:3.6}
\left[\frac{{\bf 1}^TA_{c_\alpha}^{2+k}{\bf 1}}{2{c_\alpha^2}},\frac{{\bf 1}^TA_{c_\alpha}^{4+k}{\bf 1}}{2{c_\alpha^2}},\ldots,\frac{{\bf 1}^TA_{c_\alpha}^{n-1+k}{\bf 1}}{2{c_\alpha^2}}\right]{\bf u}\equiv0\pmod{2}.
\]
Define
\begin{align*}
    M_3:=\left[
          \begin{array}{ccccc}
            \frac{{\bf 1}^TA_{c_\alpha}^2{\bf 1}}{2{c_\alpha}} & \frac{{\bf 1}^TA_{c_\alpha}^4{\bf 1}}{2{c_\alpha}} & \frac{{\bf 1}^TA_{c_\alpha}^6{\bf 1}}{2{c_\alpha}} & \cdots & \frac{{\bf 1}^TA_{c_\alpha}^{n-1}{\bf 1}}{2{c_\alpha}} \\
            \frac{{\bf 1}^TA_{c_\alpha}^3{\bf 1}}{2{c_\alpha^2}} & \frac{{\bf 1}^TA_{c_\alpha}^{5}{\bf 1}}{2{c_\alpha^2}} & \frac{{\bf 1}^TA_{c_\alpha}^{7}{\bf 1}}{2{c_\alpha^2}} & \cdots & \frac{{\bf 1}^TA_{c_\alpha}^{n}{\bf 1}}{2{c_\alpha^2}} \\
            \frac{{\bf 1}^TA_{c_\alpha}^4{\bf 1}}{2{c_\alpha^2}} & \frac{{\bf 1}^TA_{c_\alpha}^{6}{\bf 1}}{2{c_\alpha^2}} & \frac{{\bf 1}^TA_{c_\alpha}^{8}{\bf 1}}{2{c_\alpha^2}} & \cdots & \frac{{\bf 1}^TA_{c_\alpha}^{n+1}{\bf 1}}{2{c_\alpha^2}} \\
            \vdots & \vdots & \vdots & \ddots & \vdots \\
            \frac{{\bf 1}^TA_{c_\alpha}^{n+1}{\bf 1}}{2{c_\alpha^2}} & \frac{{\bf 1}^TA_{c_\alpha}^{n+3}{\bf 1}}{2{c_\alpha^2}} & \frac{{\bf 1}^TA_{c_\alpha}^{n+5}{\bf 1}}{2{c_\alpha^2}} & \cdots & \frac{{\bf 1}^TA_{c_\alpha}^{2n-2}{\bf 1}}{2{c_\alpha^2}} \\
          \end{array}
        \right]
        =\frac{\tilde{W}_{{\alpha}}^T\hat{W}'_{{\alpha}}}{2}.
\end{align*}
In view of \eqref{eq:3.6}, one has $M_3{\bf u}\equiv{\bf 0}\pmod{2}.$ By Lemma \ref{lem3.7}, we obtain that $\rank_2\left(\frac{\tilde{W}_{{\alpha}}^T\hat{W}'_{{\alpha}}}{2}\right)=\frac{n-1}{2},$ that is, the column rank of $M_3$ is full over $\mathbb{F}_2.$ Hence, ${\bf u}\equiv{\bf 0}\pmod{2},$ a contradiction.

Combining Cases 1-2, one obtains that $l$ is odd. This completes the proof.
\end{proof}

\section{\normalsize Examples}\setcounter{equation}{0}
In this section, we give two specific examples of DGA$_\alpha$S graphs which satisfies the conditions of Theorem~\ref{thm1.1} for some $\alpha\in[0,1)$.
\begin{ex}
{\rm Let the adjacency matrix of graph $G$ be given as follows:
\begin{equation*}
    A(G)=\left(
           \begin{array}{cccccccccccccc}
             0 & 1 & 0 & 0 & 0 & 0 & 1 & 0 & 1 & 0 & 0 & 0 & 0 & 1 \\
             1 & 0 & 1 & 0 & 0 & 1 & 0 & 1 & 1 & 1 & 0 & 1 & 1 & 0 \\
             0 & 1 & 0 & 0 & 0 & 1 & 0 & 0 & 0 & 0 & 0 & 1 & 0 & 0 \\
             0 & 0 & 0 & 0 & 0 & 0 & 0 & 0 & 0 & 1 & 1 & 0 & 1 & 0 \\
             0 & 0 & 0 & 0 & 0 & 1 & 0 & 1 & 1 & 1 & 0 & 0 & 0 & 1 \\
             0 & 1 & 1 & 0 & 1 & 0 & 0 & 1 & 1 & 0 & 0 & 1 & 1 & 0 \\
             1 & 0 & 0 & 0 & 0 & 0 & 0 & 1 & 1 & 0 & 1 & 0 & 0 & 1 \\
             0 & 1 & 0 & 0 & 1 & 1 & 1 & 0 & 1 & 1 & 1 & 1 & 1 & 0 \\
             1 & 1 & 0 & 0 & 1 & 1 & 1 & 1 & 0 & 0 & 1 & 0 & 0 & 0 \\
             0 & 1 & 0 & 1 & 1 & 0 & 0 & 1 & 0 & 0 & 0 & 0 & 1 & 0 \\
             0 & 0 & 0 & 1 & 0 & 0 & 1 & 1 & 1 & 0 & 0 & 1 & 1 & 0 \\
             0 & 1 & 1 & 0 & 0 & 1 & 0 & 1 & 0 & 0 & 1 & 0 & 1 & 1 \\
             0 & 1 & 0 & 1 & 0 & 1 & 0 & 1 & 0 & 1 & 1 & 1 & 0 & 0 \\
             1 & 0 & 0 & 0 & 1 & 0 & 1 & 0 & 0 & 0 & 0 & 1 & 0 & 0 \\
           \end{array}
         \right)_{14\times14}.
\end{equation*}
Hence, $D(G)={\rm diag}(4, 8, 3, 3, 5, 7, 5, 9, 7, 5, 6, 7, 7, 4).$ 
It can be computed directly using \textit{Mathematica} \cite{math} that
$$
  \det(\tilde{W}_{\frac{3}{4}}(G))=2^{7}\times5\times331\times143807\times545912603\times30283875584713\times778268539694081846899
$$
and 
$$
  \det(\tilde{W}_{\frac{5}{6}}(G))={2^{7}}\times13\times31\times37\times327773499972443320387744582054393134299875049186710656493725761.
$$
By Theorem \ref{thm1.1}, we know that $G$ is DGA$_{\frac{3}{4}}$S and DGA$_{\frac{5}{6}}$S.}
\end{ex}
\begin{ex}
{\rm Let the adjacency matrix of graph $G$ be given as follows:
\begin{equation*}
    A(G)=\left(
           \begin{array}{ccccccccccccc}
             0 & 1 & 0 & 0 & 1 & 0 & 1 & 1 & 1 & 0 & 1 & 1 & 0 \\
             1 & 0 & 0 & 0 & 1 & 0 & 0 & 0 & 0 & 0 & 1 & 0 & 0 \\
             0 & 0 & 0 & 0 & 0 & 0 & 0 & 0 & 1 & 1 & 0 & 1 & 0 \\
             0 & 0 & 0 & 0 & 1 & 0 & 1 & 1 & 1 & 0 & 0 & 0 & 1 \\
             1 & 1 & 0 & 1 & 0 & 0 & 1 & 1 & 0 & 0 & 1 & 1 & 0 \\
             0 & 0 & 0 & 0 & 0 & 0 & 1 & 1 & 0 & 1 & 0 & 0 & 1 \\
             1 & 0 & 0 & 1 & 1 & 1 & 0 & 1 & 1 & 1 & 1 & 1 & 0 \\
             1 & 0 & 0 & 1 & 1 & 1 & 1 & 0 & 0 & 1 & 0 & 0 & 0 \\
             1 & 0 & 1 & 1 & 0 & 0 & 1 & 0 & 0 & 0 & 0 & 1 & 0 \\
             0 & 0 & 1 & 0 & 0 & 1 & 1 & 1 & 0 & 0 & 1 & 1 & 0 \\
             1 & 1 & 0 & 0 & 1 & 0 & 1 & 0 & 0 & 1 & 0 & 1 & 1 \\
             1 & 0 & 1 & 0 & 1 & 0 & 1 & 0 & 1 & 1 & 1 & 0 & 0 \\
             0 & 0 & 0 & 1 & 0 & 1 & 0 & 0 & 0 & 0 & 1 & 0 & 0 \\
           \end{array}
         \right)_{13\times13}.
\end{equation*}
Clearly, $D(G)={\rm diag}(7, 3, 3, 5, 7, 4, 9, 6, 5, 6, 7, 7, 3).$ 
It can be computed directly using \textit{Mathematica} \cite{math} that
$$
  \det(\tilde{W}_{\frac{2}{3}}(G))={2^{6}}\times5\times97\times1367\times10067\times118189\times132430201\times145112609,
$$
and
\begin{align*}
  \det(\tilde{W}_{\frac{10}{11}}(G))=&{2^{6}}\times3\times2567\times3251\times18593\times110574553\\
  &\times19912837250380292202346041446775471026303813.
\end{align*}
By Theorem \ref{thm1.1}, we know that $G$ is DGA$_{\frac{2}{3}}$S and DGA$_{\frac{10}{11}}$S.}
\end{ex}
\section{\normalsize Concluding remarks}
In this paper, we give a simple arithmetic criterion for determining whether a graph $G$ is DGA$_\alpha$S. 
Obviously, the main results in \cite{0005} and \cite{0002} (i.e., Corollary \ref{cor1.01} and Corollary \ref{cor1.02}) are the direct consequences of our results in this paper.

Notice that we do not consider the case that $n$ is even and $c_\alpha(\geqslant 3)$ is odd in Theorem \ref{thm1.1}. Now, we pose the following conjecture.
\begin{conj}\label{conj1}
Let $G$ be in $\mathfrak{F}_n$ and let $\alpha$ be in $[0,1)$ with even $n$ and odd $c_\alpha\,(\geqslant 3)$. Then $G$ is DGA$_{\alpha}$S.
\end{conj}
In view of Lemma~\ref{thm2.2}, we know that in order to prove this conjecture, it suffices to show that the condition in Conjecture \ref{conj1} implies $l(U)=1$ for all matrices $U\in \Gamma_\alpha(G),$  which is equivalent to prove that any prime $p$ is not a divisor of $l(U).$ Based on Theorem \ref{thm3.1}, one obtains that any old prime $p\nmid l(U)$ for all matrices $U\in \Gamma_\alpha(G).$ Hence, in order to prove Conjecture \ref{conj1}, it is sufficient to show that $l(U)$ is odd for each  matrix $U\in \Gamma_\alpha(G).$

As it turns out in our paper, the arithmetic properties of $\det \tilde{W}_\alpha(G)$ are closely related to whether a given graph is DGA$_\alpha$S for $\alpha\in[0,1)$. By comparing Theorem~\ref{thm1.1} with Corollaries~\ref{cor1.01}-\ref{cor1.02}, we pose the following open problem:
\begin{pb}
Let $G$ be a graph with order $n$ and let $\alpha\in[0,1).$ If $\frac{\det \tilde{W}_\alpha(G)}{2^{\lfloor\frac{n}{2}\rfloor}}$ is an odd and square-free integer, determining whether $G$ is DGA$_\alpha$S or not?
\end{pb}

\end{document}